\documentclass[pdflatex,sn-mathphys-num]{sn-jnl}

\usepackage{graphicx}%
\usepackage{multirow}%
\usepackage{amsmath,amssymb,amsfonts}%
\usepackage{amsthm}%
\usepackage{mathrsfs}%
\usepackage[title]{appendix}%
\usepackage{xcolor}%
\usepackage{textcomp}%
\usepackage{manyfoot}%
\usepackage{booktabs}%
\usepackage{algorithm}%
\usepackage{algorithmicx}%
\usepackage{algpseudocode}%
\usepackage{listings}%
\usepackage{tikz}
\usetikzlibrary{shapes.geometric}
\usepackage{xcolor} 

\lstdefinestyle{pythonstyle}{
    language=Python,
    basicstyle=\ttfamily\footnotesize,
    keywordstyle=\bfseries\color{blue},
    stringstyle=\color{purple},
    commentstyle=\color{gray},
    numbers=none,
    numberstyle=\tiny\color{gray},
    stepnumber=1,
    numbersep=10pt,
    showstringspaces=false,
    frame=single,
    breaklines=true
}

\lstdefinestyle{conjecturestyle}{
    basicstyle=\ttfamily\footnotesize,
    keywordstyle=\bfseries\color{blue},
    commentstyle=\color{gray},
    showstringspaces=false,
    frame=single,
    breaklines=true,
    numbers=none,
    numberstyle=\tiny\color{gray},
    xleftmargin=0.05\textwidth, 
    xrightmargin=0.05\textwidth
}

\newtheorem{thm}{Theorem}
\newtheorem{corollary}{Corollary}

\newtheorem{conj}{Conjecture}

\newcommand{\NP}{\mathcal{NP}}
\DeclareMathOperator{\Sharp}{Sharp}

\begin{document}

\title[Article Title]{Automated Conjecturing with \emph{TxGraffiti}}

\author[1, 2]{\fnm{Randy} \sur{Davila}}\email{rrd6@rice.com}


\affil[1]{\orgname{First Principles}, \orgaddress{\street{100 University Ave}, \city{Toronto}, \state{ON M5J 1V6}, \country{Canada}}}

\affil[2]{\orgdiv{Department of Computational Applied Mathematics \& Operations Research}, \orgname{Rice University}, \orgaddress{\street{6100 Main Street}, \city{Houston}, \postcode{77005}, \state{Texas}, \country{USA}}}

\abstract{
\emph{Texas Graffiti} (\emph{TxGraffiti}) is an automated conjecturing system that operates on a finite, versioned \emph{snapshot table} of mathematical objects equipped with precomputed numerical invariants and Boolean predicates.  For a chosen target invariant and hypothesis predicate, \emph{TxGraffiti} searches over predicate--predictor combinations and generates \emph{table-true} conditional inequalities by fitting coefficients within simple, human-readable templates (e.g., univariate affine upper and lower bounds) via a sequence of small linear programs.  The resulting candidates are ranked and compressed using heuristics that measure informativeness (e.g., touch/sharp behavior) and remove redundancy or transitive implication, producing a compact set of interpretable conjectures suitable for further testing and proof. We describe the system architecture, data curation workflow, optimization models, and filtering procedures, and we provide an interactive web interface together with an open-source Python implementation.  Although our examples focus on graph invariants, the snapshot-table paradigm and optimization-and-filtering pipeline apply broadly whenever objects admit computable features and predicates.
}

\date{}
\maketitle


\section{Introduction}
\label{sec:introduction}

In 1948, Turing suggested that digital computers could play a significant role in intellectual work—including parts of mathematical research—especially in settings requiring substantial reasoning but little interaction with the external world~\cite{Turing}. This vision helped motivate early work in computer-assisted mathematics. One of the first notable efforts was Newell and Simon’s \emph{Logic Theorist}, developed in the mid-1950s~\cite{SimonNewell}. By proving many theorems from \emph{Principia Mathematica}, \emph{Logic Theorist} helped establish a template for later automated reasoning systems, even if its early achievements were modest by modern standards. Newell and Simon also advanced an influential long-range prediction that computers would eventually discover and prove important mathematical theorems~\cite{SimonNewell}. Along with many programs of its era, \emph{Logic Theorist} focused primarily on automated theorem proving, a domain that reached a widely recognized milestone with McCune’s computer-assisted proof of the Robbins conjecture in the mid-1990s~\cite{Robbins}.

While automated theorem proving matured into a well-established area, the complementary problem of computer-assisted conjecture generation (\emph{automated conjecturing}) was already being explored, though it attracted far less sustained attention. In the late 1950s, Hao Wang’s \emph{Program II} generated a large number of formulas in propositional logic—outputs that can be viewed as conjectures or candidate theorems~\cite{WangProgramII}. A central challenge for early conjecture-making systems, including Wang’s, was sifting through the resulting flood of statements to identify those of genuine mathematical interest.

A particularly influential response to this “too-many-statements” problem was Fajtlowicz’s \emph{Graffiti}, which generated mathematical conjectures by proposing inequalities among real-valued invariants of mathematical objects (see Fajtlowicz’s series \emph{On Conjectures of Graffiti}~\cite{Fajtlowicz-DM-1988, Fajtlowicz-1987, Fajtlowicz-III-1988, Fajtlowicz-IV-1990, Fajtlowicz-V-1995}). To control redundancy and avoid conjectures that are subsumed by stronger ones, Fajtlowicz introduced the \emph{Dalmatian heuristic}, a selection rule that retains a new conjecture only when it yields genuinely new information relative to those already accepted~\cite{Larson}. The program was named \emph{Graffiti} because its conjectures were metaphorically “written on the wall” for mathematicians to explore~\cite{GraffitiD}. \emph{Graffiti}’s conjectures have inspired a substantial body of follow-up work in graph theory and related areas (see DeLaViña's historical account of \emph{Graffiti}~\cite{GraffitiD}).

This paper describes \emph{Texas Graffiti} (\emph{TxGraffiti}), a heuristic- and optimization-based automated conjecturing program developed and maintained by the author since 2016 and available as an interactive web application\footnote{\url{https://txgraffiti.streamlit.app}} and PyPI-distributed Python package \texttt{txgraffiti}~\cite{txgraffiti_pypi}. The name honors \emph{Graffiti} and its successor \emph{Graffiti.pc}~\cite{Graffitipc}. Like its predecessors, \emph{TxGraffiti} primarily generates conjectures in inequality form—upper and lower bounds on a chosen target invariant under explicit Boolean hypotheses—and when a table-true upper bound and a table-true lower bound coincide on the available data, it reports the induced equality conjecture. These statements are conjectures in the traditional mathematical sense: they are validated against a finite \emph{snapshot table} of instances and may fail on objects not yet represented. Accordingly, throughout the paper we use ``table-true'' (or ``true on the snapshot'') to mean ``no counterexample among the stored instances,'' and we reserve ``theorem'' exclusively for statements established by a formal proof.

The remainder of the paper is structured as follows. In
Section~\ref{sec:related-work}, we recall the original \emph{Graffiti} program and describe its process. In Section~\ref{sec:others}, we briefly discuss other
notable conjecturing systems. In Section~\ref{sec:methods}, we detail the implementation of \emph{TxGraffiti}. In Section~\ref{sec:workflow-total-dom}, we illustrate a typical workflow to interpret and refine the conjectures produced by \emph{TxGraffiti}. Section~\ref{sec:results} presents conjectures generated by \emph{TxGraffiti} that have led to publications. Finally, Section~\ref{sec:conclusion} offers concluding remarks.

\section{Fajtlowicz's \emph{Graffiti}}
\label{sec:related-work}

Fajtlowicz’s \emph{Graffiti} is one of the first and most influential automated conjecturing systems, demonstrating that a computer can generate mathematically useful conjectures while using explicit rules to suppress redundancy and trivial output. Developed in the mid-1980s, \emph{Graffiti} searches for relations—primarily inequalities—among numerical invariants evaluated on a finite database of objects (most often graphs), and records conjectures that survive both empirical testing and a built-in informativeness screen~\cite{Fajtlowicz-DM-1988, GraffitiD}. Since its inception, \emph{Graffiti}’s conjectures have motivated substantial follow-up work, especially in \emph{graph theory} (see~\cite{Fajtlowicz-Waller-1986, Chung-1988, Alon-Seymour-1989, Fajtlowicz-Clemson-1989, Beezer-1989, Favaron-1990, Favaron-Research-1991, Favaron-1991, Favaron-DM-1993, Griggs-Kleitman-1994, Fajtlowicz-McColgan-1995, Wang-1997, Firby-1997, Bollobas-Erdos-1998, Bollobas-Riordan-1998, Caro-1998, Dankelmann-Swart-Oellermann-1998, Jelen-1999, Codenotti-2000, Beezer-Riegsecker-Smith-2001, Favaron-2003, Zhang-2004, Dankelmann-Dlamini-Swart-2005, Hansen-2009, Cygan-2012, Yue-2020}) and \emph{mathematical chemistry} (see~\cite{Fowler-1997, Fowler-1998, Fowler-1999, Fajtlowicz-Larson-2003, Fajtlowicz-Buckminsterfullerene-2005, Fajtlowicz-John-Sachs-2005, Doslic-Reti-2011}). Historical accounts also record interactions with other researchers and feedback (including Erd\H{o}s) in response to early \emph{Graffiti} conjectures~\cite{Fajtlowicz-DM-1988, GraffitiD}.

Beyond its mathematical contributions, \emph{Graffiti} has figured prominently in philosophical discussions about the nature of machine creativity. Penrose famously argued in \emph{The Emperor’s New Mind} and \emph{Shadows of the Mind} that digital computers cannot possess genuine mathematical understanding~\cite{Penrose1989Emperor, Penrose1994Shadows}. In contrast, Fajtlowicz explicitly designed \emph{Graffiti} as an empirical test, echoing Hilary Putnam’s position that whether machines can discover mathematics is not a theoretical question, but an empirical one~\cite{Fajtlowicz-1987}. \emph{Graffiti}'s ability to produce thousands of plausible and sometimes deep conjectures, later explored and proved by human mathematicians, lends weight to this empirical perspective. As Horgan noted in \emph{Scientific American}, this contrast—between Penrose’s skepticism and \emph{Graffiti}’s output—epitomizes a larger philosophical divide in understanding the role of AI in discovery~\cite{Horgan1993DeathOfProof}.

What was once a largely abstract debate has taken on renewed urgency with the emergence of large language models (LLMs) and generative AI. These systems increasingly exhibit capabilities once thought exclusive to human reasoning—such as proof synthesis inference-making the notion of fully automated machine-driven discovery more plausible than ever. Although questions remain--particularly about whether such systems truly ``understand''—automated conjecturing systems provide structured evidence that machines can participate meaningfully in open-ended mathematical inquiry. 

At a high level, \emph{Graffiti} operates on (i) a finite collection of objects (typically graphs) and (ii) a library of numerical invariants evaluated on those objects. Given a user-selected \emph{target} invariant (denoted abstractly by $\alpha$), the system searches for inequalities of the form
\[
\alpha \le f(\text{invariants}) \qquad \text{or} \qquad \alpha \ge f(\text{invariants}),
\]
where $f$ is a low-complexity expression built from the available invariants and a fixed menu of algebraic operations. Candidate expressions are generated within a grammar/complexity budget and then tested against the stored instances; candidates with a counterexample in the current database are discarded.

The key mechanism that controls output is the \emph{Dalmatian heuristic}, introduced by Fajtlowicz to avoid being overwhelmed by the enormous number of database-true but uninformative statements that arise from systematic enumeration~\cite{Larson}. In its basic form, it combines two ideas. First, a \emph{truth} screen retains only candidates with no counterexample in the current database; equivalently, using our terminology, the inequality is \emph{table-true} in the database. Second, a \emph{significance} (informativeness) screen suppresses candidates that add no new information relative to the current conjecture bank (e.g., candidates that are everywhere weaker than what is already retained on the finite database). Historical descriptions emphasize precisely this ``non-informativeness'' rejection rule together with the maintenance step that revises the conjecture bank when a new statement is accepted (including dropping conjectures that become redundant)~\cite{Fajtlowicz-V-1995, GraffitiD}. For a modern exposition and refinements, we refer the reader to~\cite{Larson}.

A convenient formalization of the significance intuition uses \emph{touch} (tightness) in the database. Fix a target invariant $\alpha$ and suppose that we are conjecturing lower bounds. Let $\mathcal{O}$ be the current database and $\mathcal{F}=\{f_1,\dots, f_k\}$ the accepted lower bound functions in the conjecture bank, where each $f_i$ satisfies
\begin{equation}\label{eq:dalmatian-lb}
\alpha(G)\ \ge\ f_i(G)\qquad\text{for all }G\in\mathcal{O}.
\end{equation}
The bank induces the aggregate bound
\begin{equation}\label{eq:aggregate-max}
B_{\mathcal{F}}(G)\ :=\ \max_{1\le i\le k} f_i(G).
\end{equation}
A graph $G\in\mathcal{O}$ is \emph{touched} if $\alpha(G)=B_{\mathcal{F}}(G)$,
or equivalently, if $\alpha(G)=f_i(G)$ for some $i$. A natural database-dependent notion of ``coverage'' is therefore:
\begin{equation}\label{eq:coverage-halting}
\bigcup_{i=1}^k \{G\in\mathcal{O}:\alpha(G)=f_i(G)\}\ =\ \mathcal{O}.
\end{equation}
In this language, the significance screen can be viewed as preferring conjectures that increase coverage (by touching at least one previously untouched object) and discarding those that never improve the bank on the finite database. Absent an explicit coverage goal such as~\eqref{eq:coverage-halting}, a generate--test-retain loop has no canonical termination point: as the expression grammar expands, one can often continue to produce new table-true candidates that are not strictly dominated by the current bank. Consequently, practical implementations enforce external stopping rules (e.g., time limits, complexity caps, or explicit coverage objectives).

\subsection{\emph{Graffiti.pc}, a variant of \emph{Graffiti}}
\label{subsec:graffitipc}

\emph{Graffiti.pc}, developed by DeLaVi\~na, is an important descendant of Fajtlowicz’s \emph{Graffiti} that preserves the same overall paradigm: given a finite database of graphs with precomputed invariant values, it searches for conjectured inequalities among algebraic expressions in those invariants and filters the resulting stream to retain a compact, nonredundant set of ``interesting'' statements~\cite{Graffitipc}. DeLaVi\~na describes a three-part architecture consisting of (i) \texttt{BuildDbs}, which computes and stores invariant values; (ii) \texttt{dalmatians}, which implement a Graffiti-style selection rule; and (iii) a graphical user interface that supports interactive choices of invariants, relations, and parameters, and also supports incorporating counterexamples provided by the user~\cite{Graffitipc}.

A distinguishing design feature is that conjecture candidates are built inside a user-controlled \emph{$\Sigma$-algebra} of expressions: invariant symbols are treated as atomic terms, and new expressions are formed by applying operations from a fixed (user-selectable) menu and representing the result as a syntax tree~\cite{Graffitipc}. In practice, the search space is regulated by syntactic constraints such as which operations are enabled, the maximum expression size/depth, and which invariants are allowed to appear. The conjectures are then generated by testing the relations (such as $\le$ and $\ge$) between the generated expressions against the stored database values~\cite{Graffitipc}. This ``grammar-first'' approach differs from statistical fitting: candidates are produced by explicit enumeration in a structured hypothesis space and then validated (or refuted) directly on the database.

As in \emph{Graffiti}, the \emph{Dalmatian} component provides the main nonredundancy control.  For inequality mining, it retains a conjecture only if it is not dominated by the current bank in the database-dependent sense that it contributes new tightness information---e.g., it touches at least one graph not already explained as sharp by the bank~\cite{Graffitipc}. Because the candidate space is, in principle, unbounded under iterative expression-building, \emph{Graffiti.pc} relies on external controls (limits on expression generation) together with database-dependent objectives; a common and conceptually natural objective is to continue until the retained conjectures collectively touch all (or nearly all) stored graphs, since candidates found beyond that point tend to be redundant from the ``explaining sharpness'' viewpoint~\cite{Graffitipc}.  Operationally, this also explains why runs can be long (potentially hours or longer) when generous operation menus and size bounds are enabled, since the number of syntactically distinct expressions -- and therefore candidate inequalities to test -- grows rapidly.

In addition to inequality mining, \emph{Graffiti.pc} supports conjecturing \emph{sufficient conditions} for graph properties via the \emph{Sophie} heuristic. Here, the output is an implication of the form
\[
\text{(simple condition)} \ \Rightarrow\ \text{(property)},
\]
validated against the current database and prioritized toward conditions that are broadly correct on known instances while remaining syntactically simple.  In early \emph{Graffiti.pc} case studies, Sophie-generated conjectures were used, for example, to propose structural descriptions of equality cases and related sufficient-condition phenomena arising in total domination~\cite{totaldom2007}.

\subsection{\emph{TxGraffiti} in context: comparison with \emph{Graffiti} and \emph{Graffiti.pc}}
\label{subsec:txgraffiti-vs-graffiti}

\emph{TxGraffiti} descends from \emph{Graffiti} and \emph{Graffiti.pc} in the sense that it targets the same style of output: concise conjectured relationships among invariants, together with explicit preferences for nonredundancy and mathematical usefulness. However, in this paper, \emph{TxGraffiti} refers to an algorithmic design to convert a finite collection of objects—equipped with numerical invariants and Boolean predicates—into a filtered list of candidate inequalities. The interactive website is one system-level instantiation of this design for graph invariants, while the companion Python package \texttt{txgraffiti} factors the design into modular components so it can be applied to arbitrary snapshot tables~\cite{txgraffiti_pypi}.

The first difference concerns how evidence is represented and reused.  Both \emph{Graffiti} and \emph{Graffiti.pc} test candidate statements against an evolving database of objects for which invariant values are computed (and stored) by program-specific workflows.  In \emph{TxGraffiti}, on the contrary, the evidence is made explicit as a \emph{snapshot table} $T$: rows are objects and columns are precomputed invariants and predicates.  Conjecturing then reduces to a sequence of structured queries and optimization calls over $T$, and the table itself becomes a reusable artifact that can be inspected, versioned, and shared.

The second difference concerns the generation of candidates. \emph{Graffiti} and \emph{Graffiti.pc} search within hand-designed or user-configured expression grammars (\emph{expression trees}), proposing many candidate formulas and testing them against stored instances. In \emph{TxGraffiti}, the algorithm instead fixes a simple parametric hypothesis class (for example, affine inequalities of the form $\alpha \le m x + b$ under a predicate $P$) and fits the coefficients by solving an optimization model over the corresponding rows of the snapshot table.  For a chosen target invariant $\alpha$, it enumerates predictor--predicate combinations available in the current schema and solves the associated models.  

A third difference concerns \emph{termination} and the role of Dalmatian-style filtering. In the original \emph{Graffiti} paradigm, the Dalmatian heuristic functions as an internal selection rule within an open-ended search loop: absent an explicit external budget, target list, or coverage objective, there is no canonical point at which the pipeline must stop, since one can often continue to generate additional table-true statements on a finite database that are not strictly dominated by the current bank. \emph{TxGraffiti} also uses Dalmatian-style nonredundancy filters (see Section~\ref{sec:static-dalmatian}), but its conjecturing runs terminate by construction: for a fixed snapshot table $T$, a fixed predicate vocabulary, and a fixed hypothesis class (together with an explicit complexity budget), the system solves a finite batch of optimization problems and returns the resulting table-true candidates after redundancy filtering.

Direct head-to-head benchmarking against historical installations of \emph{Graffiti} and \emph{Graffiti.pc} is difficult to interpret, since reported behavior depends strongly on the underlying databases, invariant libraries, operation menus, and interactive workflows.  In addition, the systems are organized around different search regimes: \emph{TxGraffiti} runs a finite batch of optimization-and-filtering steps over a fixed snapshot table and template family, whereas \emph{Graffiti}/\emph{Graffiti.pc} were designed for open-ended conjecture banking under user-driven exploration.  With these caveats in mind, informal comparisons indicate that \emph{TxGraffiti} reproduces many classical Graffiti-style conjectures whenever invariant vocabularies and hypothesis templates overlap, while also producing qualitatively different conjectures when optimization-based templates and richer predicate vocabularies are enabled. For example, using invariants available in earlier systems, \emph{TxGraffiti} discovered the bound that the independence number is at most the matching number for $r$-regular graphs (see~\cite{CaDaPe2020}).

\section{Other Notable Conjecturing Systems}
\label{sec:others}

Automated conjecturing systems have developed along several largely independent lines, ranging from early AI-driven concept formation, to graph-focused optimization and inequality mining, and (more recently) statistical and
machine-learning pipelines. We briefly survey representative systems in this broader landscape.

\subsection{Symbolic systems}

Early automated conjecturing systems largely pursued symbolic concept formation and rule-based exploration, often with explicit measures of ``interestingness.'' A later and particularly influential branch focused on
discrete structures---most notably graphs---where conjectures can be generated from computed invariants and tested systematically against growing collections of examples and counterexamples. We briefly highlight representative systems in this lineage.

Lenat's \emph{AM}~\cite{Lenat_1, Lenat_2, Lenat_3} (1970s) explored elementary mathematical concept formation via heuristic search. Starting from a small set of primitives and production rules, \emph{AM} generated new concepts and conjectural relationships among them, using explicit heuristics and measures of ``interestingness'' to prioritize which avenues to pursue~\cite{Lenat_2, Lenat_1}. Although its mathematical output was limited in scope, \emph{AM} helped establish the general paradigm of automated discovery as a guided expansion of a concept space under heuristic control.

Epstein's \emph{GT} (\emph{Graph Theorist})~\cite{Epstein_2} (1980s) focused directly on graph theory. From the definitional information in its knowledge base, \emph{GT} constructs examples, proposes conjectures, and (in some cases) proves theorems about graph-theoretic concepts~\cite{Epstein_2}. It is an early example of an AI system organized around a specific mathematical domain and illustrates how example generation and property evaluation can be coupled to conjecture formation.

Colton's \emph{HR}~\cite{Colton_1, Colton_2, Colton_3} (late 1990s) generalized the concept-formation paradigm across multiple domains, including algebra, graph theory, and number theory. Given foundational domain information supplied by the user (e.g., axioms, objects, or tables), \emph{HR} applies production rules to construct new concepts and then proposes conjectures by detecting regularities and shared structure between independently generated concepts, guided by explicit measures of ``interestingness''~\cite{Colton_1, Colton_3}.

Hansen and Caporossi’s \emph{AutoGraphiX} (\emph{AGX})~\cite{AGX_1, AGX_2} (developed since the late 1990s) emphasizes combinatorial optimization on graph spaces. It frames the search for extremal graphs (for a chosen invariant expression under constraints) as an optimization problem and applies Variable Neighborhood Search, using neighborhoods defined by edge additions, deletions, exchanges, and related local transformations to escape local optima and identify extremal candidates~\cite{AGX_1}. This workflow is particularly effective in extremal graph theory, where sharpness phenomena and small families of extremal examples often drive conjecture formation. In addition to producing extremal examples, AGX’s “automation” mode supports conjecture discovery by repeatedly solving extremal instances across related objectives, constraint choices, and recording patterns that persist between runs~\cite{AGX_2}.

M\'elot’s \emph{GraPHedron}~\cite{graphedron_1} (2008) introduced a geometric perspective: represent each graph as a point in a space of invariants and study the convex hull of these points. Facets of the resulting polytope correspond to linear inequalities among invariants, and extremal graphs appear as vertices supporting tight inequalities~\cite{graphedron_1}. This polyhedral view supports the extraction of compact, nonredundant inequality libraries and has been extended in successor projects such as \emph{PHOEG}, which explicitly builds on and augments the GraPHedron approach~\cite{devillez2019}.

Larson and Van Cleemput’s \emph{Conjecturing} program~\cite{Larson, LarsonVC2017} (mid-2010s), implemented within \emph{SageMath}, revisits automated conjecturing in a flexible user-facing environment. In their inequality-mining work, the program searches for conjectured relations between real-valued invariants (often under hypotheses) and filters candidates using Graffiti-style screening, implemented via a truth test together with a Dalmatian-style significance test~\cite{Larson}. Complementing this, they also study \emph{property-relation conjectures}—relations among properties of objects—aiming to output simple statements that hold for all objects currently known to the program, with case studies including conjectures about \emph{Hamiltonicity}~\cite{LarsonVC2017}. The same conjecturing architecture has also been demonstrated outside of graph theory, for example, in combinatorial game theory (\emph{Chomp}), where the goal is a compact and human-readable “theory” of positions and the authors prove one generated conjecture as a demonstration of output quality~\cite{BradfordEtAlChomp2020}. As they emphasize, this Dalmatian-style screening still requires an explicit stopping condition; an idealized endpoint is exact predictive power on the stored instance set, while in general one must impose an external cap such as time or complexity~\cite{BradfordEtAlChomp2020}. Related community outreach efforts include the \emph{Independence Number Project}\footnote{\url{https://independencenumber.wordpress.com/}}.

\subsection{Statistical and machine learning approaches}\label{subsec:machine_learning}

The systems above primarily rely on symbolic rules, combinatorial optimization, or geometric representations of invariant spaces. A complementary line of work uses statistical learning to detect regularities in large mathematical datasets and to propose candidate relationships that can then be examined, refined, and
(when possible) proved. In this paradigm, the immediate output is often empirical evidence or a compact candidate pattern rather than a deductive proof, but the resulting hypotheses can nevertheless play the same role as traditional conjectures: they focus attention, suggest new definitions, and guide proof search.

The \emph{Ramanujan Machine} of Raayoni et al.~\cite{Raayoni2021RamanujanMachine} is a data-driven discovery engine for conjecturing new formulas for fundamental
constants (e.g., $\pi$, $e$, Catalan's constant, and special values of $\zeta(s)$), with output prominently featuring structured continued-fraction representations. The system searches within restricted parametric families encoded by its algorithms, including a meet-in-the-middle procedure and a gradient-descent scheme tailored to the recurrent structure of continued fractions, and retains candidates by high-precision numerical agreement rather than by proof~\cite{Raayoni2021RamanujanMachine}. Accordingly, the output is conjectural in an explicitly experimental sense. Subsequent work has rigorously proved and, in many cases, generalized families of conjectures arising from this program~\cite{Yamamoto2024RamanujanMachineProof, DoughertyBlissZeilberger2023}.

DeepMind’s 2021 study of Davies et al.~\cite{Davies2021AdvancingMathAI}
illustrated a different ML-assisted discovery loop: train predictive models on large datasets of mathematical objects, then use attribution methods to isolate the salient structure that can be reformulated into human-legible conjectures and, in favorable cases, theorems~\cite{Davies2021AdvancingMathAI}. In knot theory, they trained a supervised model to predict the signature $\sigma(K)$ of hyperbolic knots from geometric features and used gradient-based attribution to identify cusp-shaped quantities—specifically meridional translation $\mu$ and
longitudinal translation $\lambda$—as dominant predictors~\cite{Davies2021AdvancingMathAI}. This guided the introduction of the \emph{natural slope} $\mathrm{Re}(\lambda/\mu)$ and led to conjectures and a proved theorem relating slope and signature, with additional geometric correction terms~\cite{Davies2021AdvancingMathAI}. In representation theory, a message-passing neural network was trained to predict Kazhdan--Lusztig polynomials from unlabeled Bruhat intervals; attribution and analysis of salient subgraphs highlighted the role of extremal reflections and suggested a hypercube-based decomposition heuristic, providing structural evidence relevant to combinatorial invariance conjecture~\cite{Davies2021AdvancingMathAI}.

In contrast to expression-grammar enumeration and to purely statistical pattern extraction, \emph{TxGraffiti} is organized around an explicit snapshot table: a finite, reusable artifact whose rows are objects and whose columns are precomputed invariants and Boolean predicates. For a fixed schema, the program generates conjectures by enumerating a finite family of optimization models indexed by target, predictor, and hypothesis predicate and then applies Graffiti-style screening to retain a compact, nonredundant set of table-true statements together with their equality witnesses. We now describe this pipeline in a concrete graph-theoretic instance.

\section{\emph{TxGraffiti} Design Principles}
\label{sec:methods}
In this section, we describe the design principles underlying \emph{TxGraffiti} and the conjecturing pipeline implemented in the current web application.\footnote{\url{https://txgraffiti.streamlit.app}} To support reproducibility, companion code for the illustrative example in this section is available in the paper’s GitHub repository.\footnote{\url{https://github.com/RandyRDavila/Automated-conjecturing-in-mathematics-with-TxGraffiti}}

For exposition, we begin with the small dataset $\mathcal{O}$ shown in Figure~\ref{fig:graphs}, consisting of nine connected graphs. These graphs are included only to make the mechanics of the pipeline and the structure of the underlying table completely explicit; they are not intended to be sufficiently rich to support high-quality conjecturing. In research use, \emph{TxGraffiti} relies on a substantially larger snapshot table (currently on the order of a few hundred) of \emph{special graphs}, curated over several years from named extremal constructions and standard counterexamples in the literature, together with graphs that arose during conjecture testing and proof attempts. The guiding principle of this curation is coverage rather than randomness: the table is strengthened by adding diagnostically important instances—especially equality (or near-equality) witnesses for candidate bounds and newly discovered counterexamples that expose missing hypotheses or structural obstructions.

This workflow differs from the standard supervised-learning setting.  The goal is not to fit a predictor and estimate generalization from held-out data, but to propose concise universal statements supported by all currently available evidence: generate a candidate conjecture, test it against every stored instance, and update the evidence base when new witnesses are found.  Nevertheless, any finite snapshot has limitations: it may underrepresent large-order behavior or omit rare obstructions, so table-truth is evidence rather than a guarantee of universal validity.  In practice, we mitigate this by repeatedly enlarging the snapshot with additional extremal families and counterexamples and rerunning the same conjecturing queries; conjectures that persist under such table growth are prioritized for follow-up, while failures are incorporated into the table and used to repair hypotheses or refine the conjecture bank.

From the illustrative dataset $\mathcal{O}$, we compute the graph-theoretic features used by \emph{TxGraffiti}, including numerical invariants (e.g., the number of vertices $n$, the matching number $\mu$, and the independence number $\alpha$) and Boolean predicates (e.g., whether a graph is bipartite or regular). Appendix~\ref{app:terminology} lists the graph-theoretic terminology and notation used throughout the paper.

\begin{figure}[htb]
\begin{center}
\begin{tikzpicture}[scale=.8,style=thick,x=1cm,y=1cm]
\def\vr{3.5pt}

\path (-5, 3) coordinate (a1);
\path (-4, 3) coordinate (a2);
\path (-3, 3) coordinate (a3);
\draw (a1) -- (a2) -- (a3);
\draw (a1) [fill=black] circle (\vr);
\draw (a2) [fill=black] circle (\vr);
\draw (a3) [fill=black] circle (\vr);
\node at (-4, 2.5) {($G_1$)};

\path (-1, 3) coordinate (b1);
\path (0, 4) coordinate (b2);
\path (1, 3) coordinate (b3);
\draw (b1) -- (b2) -- (b3) -- (b1);
\draw (b1) [fill=black] circle (\vr);
\draw (b2) [fill=black] circle (\vr);
\draw (b3) [fill=black] circle (\vr);
\node at (0, 2.5) {($G_2$)};

\path (3, 3) coordinate (c1);
\path (4, 3) coordinate (c2);
\path (4, 4) coordinate (c3);
\path (3, 4) coordinate (c4);
\draw (c1) -- (c2) -- (c3) -- (c4) -- (c1);
\draw (c1) [fill=black] circle (\vr);
\draw (c2) [fill=black] circle (\vr);
\draw (c3) [fill=black] circle (\vr);
\draw (c4) [fill=black] circle (\vr);
\node at (3.5, 2.5) {($G_3$)};

\path (-5, 0) coordinate (d1);
\path (-4, 1) coordinate (d2);
\path (-3, 0) coordinate (d3);
\path (-4, 0) coordinate (d4);
\draw (d1) -- (d2) -- (d3) -- (d1);
\draw (d2) -- (d4);
\draw (d1) [fill=black] circle (\vr);
\draw (d2) [fill=black] circle (\vr);
\draw (d3) [fill=black] circle (\vr);
\draw (d4) [fill=black] circle (\vr);
\node at (-4, -0.5) {($G_4$)};

\path (-1, 0) coordinate (e1);
\path (0, 1) coordinate (e2);
\path (1, 0) coordinate (e3);
\path (0, -1) coordinate (e4);
\draw (e1) -- (e2) -- (e3) -- (e4) -- (e1);
\draw (e1) -- (e3);
\draw (e2) -- (e4);
\draw (e1) [fill=black] circle (\vr);
\draw (e2) [fill=black] circle (\vr);
\draw (e3) [fill=black] circle (\vr);
\draw (e4) [fill=black] circle (\vr);
\node at (0, -1.5) {($G_5$)};

\path (3, 1) coordinate (f1);
\path (4, 1) coordinate (f2);
\path (5, 1) coordinate (f3);
\path (6, 1) coordinate (f4);
\path (3, -.5) coordinate (f5);
\path (4, -.5) coordinate (f6);
\path (5, -.5) coordinate (f7);
\path (6, -.5) coordinate (f8);
\draw (f1) -- (f5) -- (f2) -- (f6) -- (f3) -- (f7) -- (f4) -- (f8);
\draw (f1) -- (f6);
\draw (f1) -- (f7);
\draw (f1) -- (f8);
\draw (f2) -- (f5);
\draw (f2) -- (f7);
\draw (f2) -- (f8);
\draw (f3) -- (f5);
\draw (f3) -- (f6);
\draw (f3) -- (f8);
\draw (f4) -- (f5);
\draw (f4) -- (f6);
\draw (f4) -- (f7);
\draw (f1) [fill=black] circle (\vr);
\draw (f2) [fill=black] circle (\vr);
\draw (f3) [fill=black] circle (\vr);
\draw (f4) [fill=black] circle (\vr);
\draw (f5) [fill=black] circle (\vr);
\draw (f6) [fill=black] circle (\vr);
\draw (f7) [fill=black] circle (\vr);
\draw (f8) [fill=black] circle (\vr);
\node at (4.5, -1.5) {($G_6$)};

\path (-5, -3) coordinate (g1);
\path (-6, -4) coordinate (g2);
\path (-5, -4) coordinate (g3);
\path (-4, -4) coordinate (g4);
\draw (g1) -- (g2);
\draw (g1) -- (g3);
\draw (g1) -- (g4);
\draw (g1) [fill=black] circle (\vr);
\draw (g2) [fill=black] circle (\vr);
\draw (g3) [fill=black] circle (\vr);
\draw (g4) [fill=black] circle (\vr);
\node at (-5, -4.5) {($G_7$)};

\path (-1, -3) coordinate (h1);
\path (0, -3) coordinate (h2);
\path (-2, -4) coordinate (h11);
\path (-1, -4) coordinate (h12);
\path (0, -4) coordinate (h21);
\path (1, -4) coordinate (h22);
\draw (h1) -- (h2);
\draw (h1) -- (h11);
\draw (h1) -- (h12);
\draw (h2) -- (h21);
\draw (h2) -- (h22);
\draw (h1) [fill=black] circle (\vr);
\draw (h2) [fill=black] circle (\vr);
\draw (h11) [fill=black] circle (\vr);
\draw (h12) [fill=black] circle (\vr);
\draw (h21) [fill=black] circle (\vr);
\draw (h22) [fill=black] circle (\vr);
\node at (-.5, -4.5) {($G_8$)};

\path (3, -3) coordinate (i1);
\path (4, -4) coordinate (i2);
\path (5, -3) coordinate (i3);
\path (6, -3) coordinate (i4);
\path (7, -4) coordinate (i5);
\path (8, -3) coordinate (i6);
\draw (i1) -- (i2) -- (i3) -- (i1);
\draw (i4) -- (i5) -- (i6) -- (i4);
\draw (i3) -- (i4);
\draw (i1) [fill=black] circle (\vr);
\draw (i2) [fill=black] circle (\vr);
\draw (i3) [fill=black] circle (\vr);
\draw (i4) [fill=black] circle (\vr);
\draw (i5) [fill=black] circle (\vr);
\draw (i6) [fill=black] circle (\vr);
\node at (5.5, -4.5) {($G_9$)};

\end{tikzpicture}
\end{center}
\caption{A dataset of 9 connected graphs: ($G_1$) the path graph on 3 vertices, ($G_2$) the cycle graph on 3 vertices, ($G_3$) the cycle graph on 4 vertices, ($G_4$) the \emph{diamond graph}, ($G_5$) the complete graph $K_4$, ($G_6$) the complete bipartite graph $K_{4,4}$, ($G_7$) the star graph $K_{1,3}$, ($G_8$) the double star graph $S(2,2)$, ($G_9$) and graph obtained by two triangles joined by an edge.}
\label{fig:graphs}
\end{figure}

Table~\ref{tab:graph-data} summarizes these precomputed properties for each graph in our dataset and serves as the evidence base for \emph{TxGraffiti}'s conjecturing queries (in this toy example). This table-first design differs from \emph{Graffiti} and \emph{Graffiti.pc}, which compute (and repeatedly recompute) invariant values during candidate-generation and screening, a workflow that can lead to long runs when expensive invariants and large candidate menus are enabled. In \emph{TxGraffiti}, invariant computation is moved to an offline snapshot-construction step, and conjecturing reduces to a finite batch of optimization calls and subsequent redundancy filtering on the fixed snapshot table (see the static-Dalmatian complexity discussion in Section~\ref{sec:static-dalmatian}). Consequently, the cost of expensive invariant computation is amortized across many conjecturing queries on the same snapshot, enabling rapid interactive feedback once the table is built.\footnote{In our web deployment, computing the full snapshot table (all property columns) takes on the order of $30$--$60$ minutes, while once the snapshot is available, a typical single-target conjecturing query completes in minutes; recomputing all invariants per query would dominate interactive latency. Timing depends on hardware, the invariant set, and snapshot size.}

\begin{table}[h!]
\centering
\begin{tabular}{|c|c|c|c|c|c|c|c|c|c|c|}
\hline
\textbf{name} & \textbf{$n$} & \textbf{$\mu$} & \textbf{$\alpha$} & \textbf{$n - \mu$} & \textbf{$n - \delta$} & \textbf{$\Delta^2$} & \textbf{connected} & \textbf{tree} & \textbf{regular} & \textbf{bipartite} \\ \hline
$G_1$ & 3 & 1 & 2 & 2 & 2 & 4 & True & True & False & True \\ \hline
$G_2$ & 3 & 1 & 1 & 2 & 1 & 4 & True & False & True & False \\ \hline
$G_3$ & 4 & 2 & 2 & 2 & 2 & 4 & True & False & True & True \\ \hline
$G_4$ & 4 & 2 & 2 & 2 & 2 & 9 & True & False & False & False \\ \hline
$G_5$ & 4 & 2 & 1 & 2 & 1 & 9 & True & False & True & False \\ \hline
$G_6$ & 8 & 4 & 4 & 4 & 4 & 16 & True & False & True & True \\ \hline
$G_7$ & 4 & 1 & 3 & 3 & 3 & 9 & True & True & False & True \\ \hline
$G_8$ & 6 & 2 & 4 & 4 & 5 & 9 & True & True & False & True \\ \hline
$G_9$ & 6 & 3 & 2 & 3 & 4 & 9 & True & False & False & False \\ \hline
\end{tabular}
\caption{Precomputed numerical invariants and Boolean predicates for the illustrative dataset in Figure~\ref{fig:graphs}. Here $n=|V(G)|$ is the order, $\mu(G)$ is the matching number, $\alpha(G)$ is the independence number, $\delta(G)$ and $\Delta(G)$ are the minimum and maximum degrees, and $n-\mu$, $n-\delta$, and $\Delta^2$ are simple derived features included as additional columns.}
\label{tab:graph-data}
\end{table}

Because \emph{TxGraffiti} reports statements that are \emph{table-true} on a concrete snapshot, the dominant robustness risk is not statistical noise in the supervised-learning sense, but \emph{table corruption}: an invariant value computed incorrectly, a Boolean predicate mislabeled, or a missing value silently coerced into a numerical default.  Any such error can (by design) create spurious counterexamples or spurious sharpness events, thereby changing which candidates pass the truth test, altering touch counts, and affecting downstream pruning.  For this reason, \emph{TxGraffiti} treats the snapshot table as a versioned artifact whose entries are intended to be reproducible and auditable: invariants are computed by deterministic routines when available; expensive or heuristic computations are explicitly flagged; and conjectures are interpreted as prompts for verification, counterexample search, and table revision.  In practice, when a candidate conjecture appears sensitive to a small set of rows, those rows are inspected and recomputed, and the conjecture bank is regenerated on the corrected snapshot.  This workflow mirrors standard mathematical practice: a single counterexample is informative, but it must be a \emph{genuine} counterexample rather than a data artifact.

\subsection{Generating inequalities}
\label{subsec:generating-inequalities}

Once a snapshot table of objects with precomputed numerical invariants and Boolean predicates is available, the next step is to generate candidate inequalities that are \emph{table-true} on that snapshot.  In this subsection, we fix a target invariant—denoted abstractly by $\alpha$ (for example, the independence number $\alpha(G)$)—and, under a chosen Boolean predicate $P$, compute the affine upper and lower bounds for $\alpha$ as functions of a single predictor column $x$ (for example, $x=\mu$).  Throughout, $[N]=\{1,2,\dots,N\}$ indexes the rows of the snapshot table, and we write
\[
I(P)\ :=\ \{\,i\in[N]: P_i=\texttt{True}\,\}
\]
for the set of rows satisfying the hypothesis predicate $P$ (discarding the case $I(P)=\varnothing$).

\begin{figure}[htb]
  \centering
  \includegraphics[width=0.8\textwidth]{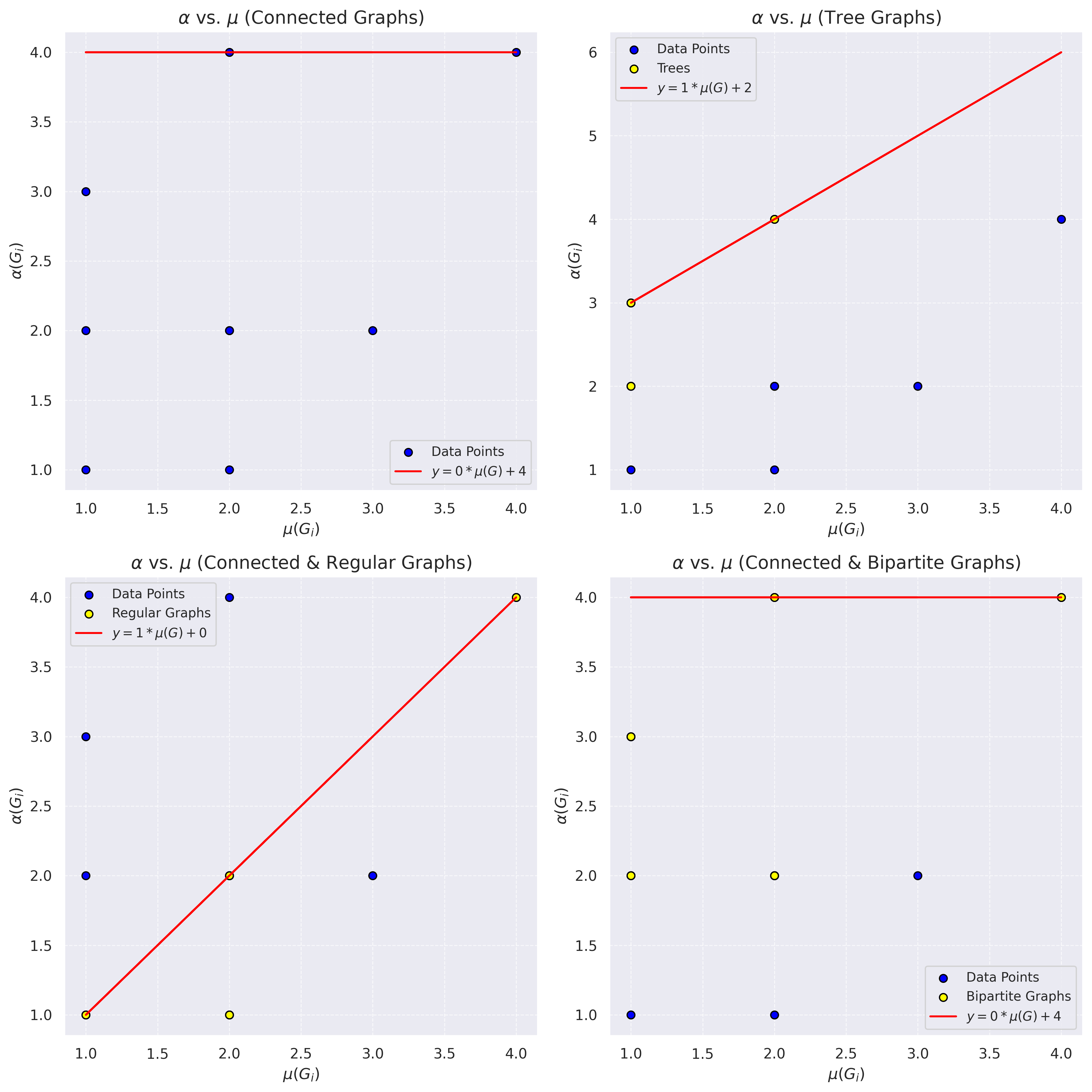}
  \caption{Independence number $\alpha$ versus matching number $\mu$ under several Boolean predicates.  Red lines show solutions to the upper-bound model~\eqref{eq:model1}.}
  \label{fig:independence_number_upper_conjectures_grid}
\end{figure}

\medskip
\noindent
To compute a table-true upper bound on $\alpha$ of the form $mx+b$ over the support set $I(P)$, we solve the linear program
\begin{equation}\label{eq:model1}
\begin{aligned}
& \underset{m, b}{\text{minimize}} && \sum_{i\in I(P)} \big(mx_i + b - \alpha_i\big) \\
& \text{subject to} && \alpha_i \le mx_i + b, \qquad \forall i \in I(P),
\end{aligned}
\end{equation}
where $x_i$ denotes the predictor value in row $i$ and $\alpha_i$ the target value.  The constraints enforce \emph{table-truth} on $I(P)$, and the objective minimizes the total \emph{slack} $mx_i+b-\alpha_i$ across the supported rows.  Thus, among all all affine functions that lie above $\alpha$ in the snapshot, \eqref{eq:model1} selects one that is, in the sense of $\ell_1$, as close as possible to $\alpha$ while remaining an upper bound.
Figure~\ref{fig:independence_number_upper_conjectures_grid} illustrates the resulting bounds for $\alpha$ versus $\mu$ in several predicates.

\medskip
\noindent
Similarly, to compute a table-true lower bound of the form $mx+b$ over the same support set, we solve
\begin{equation}\label{eq:model2}
\begin{aligned}
& \underset{m, b}{\text{maximize}} && \sum_{i\in I(P)} \big(mx_i + b - \alpha_i\big) \\
& \text{subject to} && \alpha_i \ge mx_i + b, \qquad \forall i \in I(P).
\end{aligned}
\end{equation}
Here, the constraints enforce $mx_i+b\le \alpha_i$ on $I(P)$, so every term in the objective is nonpositive; therefore, maximizing the sum chooses, among all affine lower bounds, one that stays as close as possible to $\alpha$ (again in a $\ell_1$ sense) while remaining below it.

\medskip
\noindent
For a table-true inequality restricted to $I(P)$, the \emph{touch number} is the number of rows $i\in I(P)$ on which the inequality holds with equality; these rows are the \emph{sharp instances} of the bound on the current snapshot.  In the baseline setting, we apply \eqref{eq:model1} and \eqref{eq:model2} across all choices of predictor invariant $x$ (ranging over numerical columns distinct from the target) and all predicates $P$ available in the schema.  The output of this stage is a list of table-true one-predictor affine bounds together with their touch numbers and sharp-instance sets; subsequent filtering (Section~\ref{sec:static-dalmatian}) reduces redundancy while retaining useful aliases.

\paragraph{Implementation details (LP solvers).}
All linear programs used to fit the affine bounds in \eqref{eq:model1}--\eqref{eq:model2} are modeled and solved in Python. In the companion Jupyter notebooks and in the interactive website, we use the PuLP modeling layer \cite{mitchell2011pulp} with its default backend solver CBC (COIN-OR Branch-and-Cut)~\cite{forrest2005cbc}.
In the current released Python package, the same LPs are solved using \texttt{scipy.optimize.linprog}~\cite{2020SciPy-NMeth} using the HiGHS solver~\cite{huangfu2018parallelizing}. When an LP admits multiple optimal solutions (e.g., more than one pair $(m,b)$ attains the optimum), we accept any optimal solution returned by the solver; such non-uniqueness may change the reported coefficients but does not affect the induced table-truth of the resulting conjecture, nor downstream touch/sharp sets or redundancy filtering.

\subsection{Extending single variable bounds to multivariable and nonlinear bounds}
\label{subsec:limits}

The baseline inequality generator searches, for each Boolean predicate $P$ and each numerical predictor column $x$ in the snapshot table, within a deliberately restricted but highly interpretable hypothesis class. The upper  bounds are obtained by solving~\eqref{eq:model1} (with $x_i$ set to the chosen predictor values) and the lower bounds by~\eqref{eq:model2}. Consequently, every upper-bound conjecture produced by the baseline procedure has the univariate affine form
\begin{equation}\label{eq:univariate-affine}
\alpha \le m x + b,
\end{equation}
and every lower-bound conjecture has the analogous form $\alpha \ge mx+b$, where $m,b\in\mathbb{R}$ are returned by the solver.  If $x$ is itself a derived feature (for example, $x=n-\mu$), then the right-hand side can reflect richer structure via schema-level feature engineering, although the optimization step still treats one numerical column at a time. More broadly, while the baseline generator is \emph{linear in the fitted coefficients}, it can capture many \emph{nonlinear} relationships in the underlying invariants by expanding the snapshot schema with additional derived columns (e.g., $\log n$, $x^2$, $\Delta^2$, $\min\{x,y\}$, $\max\{x,y\}$, or piecewise-linear surrogates), and then applying the same linear-programming pipeline to the enlarged feature set.  The main effect of this added expressiveness is an increase in the number of candidate predictors (and hence the number of LPs solved), trading runtime for a broader search space.  This design keeps reported conjectures easy to interpret and compare, but it excludes genuinely multivariate linear relationships unless they are included explicitly through additional columns or derived features.

\medskip
The same template extends directly to multivariate affine bounds in the newer PyPI-distributed \texttt{txgraffiti} Python package~\cite{txgraffiti_pypi}. Let $X\in\mathbb{R}^{|I(P)|\times k}$ denote the matrix of $k$ predictor columns restricted to the support set $I(P)$, and let $y\in\mathbb{R}^{|I(P)|}$ denote the corresponding target values (so $y_i=\alpha_i$).  In the implementation, multivariate bounds are obtained by a \emph{minimum total slack} linear program with explicit box constraints on coefficients and intercept, which stabilizes the fit and enforces interpretability.

For an upper bound $y \le Xw + b\mathbf{1}$, we solve
\begin{equation}\label{eq:multivar-upper-minslack}
\begin{aligned}
& \underset{w,\, b,\, s}{\text{minimize}}
& & \sum_{i\in I(P)} s_i \\
& \text{subject to}
& & X_i w + b - s_i = y_i, \qquad \forall i\in I(P),\\
& & & s_i \ge 0, \qquad \forall i\in I(P),\\
& & & -M \le w_j \le M \ \ (\forall j), \qquad -M \le b \le M,
\end{aligned}
\end{equation}
where $X_i$ is the $i$th row of $X$, $s$ is a vector of nonnegative slack variables, and $M>0$ is a coefficient bound chosen by the user. Eliminating $s_i$ shows that $X_iw+b \ge y_i$ for all $i\in I(P)$, and the objective minimizes the total slack $\sum_i (X_iw+b-y_i)$.

For a lower bound $y \ge Xw + b\mathbf{1}$, we solve the following.
\begin{equation}\label{eq:multivar-lower-minslack}
\begin{aligned}
& \underset{w,\, b,\, s}{\text{minimize}}
& & \sum_{i\in I(P)} s_i \\
& \text{subject to}
& & X_i w + b + s_i = y_i, \qquad \forall i\in I(P),\\
& & & s_i \ge 0, \qquad \forall i\in I(P),\\
& & & -M \le w_j \le M \ \ (\forall j), \qquad -M \le b \le M.
\end{aligned}
\end{equation}
Here, feasibility implies $X_iw+b \le y_i$ on $I(P)$, and minimizing $\sum_i s_i = \sum_i (y_i-(X_iw+b))$ is equivalent (for fixed $y$) to maximizing the aggregate fitted bound $\sum_i (X_iw+b)$ while remaining below $y$.

Specializing~\eqref{eq:multivar-lower-minslack} to $k=2$ predictors $(x,z)$ yields table-true inequalities of the form
\[
\alpha \ge m_1 x + m_2 z + b,
\]
capturing interactions between invariants that are inaccessible to a single predictor unless those interactions are explicitly engineered as additional derived columns.  After solving, \emph{TxGraffiti} additionally rationalizes the coefficients to small denominators and applies stricter “human-niceness” thresholds to suppress overly large or fragile expressions.

\medskip
Although the underlying software supports multivariate bounds, the interactive website defaults to single-predictor models.  In practice, many multivariate expressions were less useful to domain experts: as syntactic complexity increases, a conjecture can become harder to interpret, harder to communicate, and less likely to suggest a proof strategy, even when it is table-true and numerically tighter on a fixed snapshot.  This preference for simplicity motivates the default presentation in the web interface, which prioritizes single-predictor bounds (together with a small number of engineered features).

At the same time, the multivariate structure can remain interpretable when expressed through a natural derived target.  As an illustration, suppose that we add the derived column $\alpha \Delta$; the product of two numerical columns.
Applying the lower-bound model~\eqref{eq:multivar-lower-minslack} to $\alpha \Delta$ (as a target) over pairs of predictors yields inequalities of the form
\[
\alpha \Delta \ge m_1 a + m_2 \operatorname{res} + b,
\]
which can be rearranged into lower bounds on $\alpha$ involving a small number of invariants; in this case $a$ and $res$. For example, this procedure produced (after standard filtering) the following early (open) conjecture.
\begin{conj}[\emph{TxGraffiti}]
If $G$ is a connected graph with order $n(G) \ge 3$, then
\[
\alpha(G) \ge \frac{a(G) + \operatorname{res}(G)}{\Delta(G)}.
\]
Here, $\alpha(G)$ is the independence number, $a(G)$ is the annihilation number~\cite{Yue-2020}, $\operatorname{res}(G)$ is the Havel-Hakimi residue~\cite{Favaron-1991}, and $\Delta(G)$ is the maximum degree.
\end{conj}

This inequality illustrates the “middle ground” motivating our design: it combines a small number of standard invariants in a form that is easy to state and plausibly amenable to proof, while still expressing an interaction that cannot be captured by a single predictor column.

\subsection{The static-Dalmatian heuristic}
\label{sec:static-dalmatian}

Because \emph{TxGraffiti} systematically enumerates predictors and hypotheses—testing each numerical column as a predictor for the target under each available Boolean predicate—the baseline stage typically generates a large family of table-true candidate bounds (often dozens to hundreds for a single target invariant).  We therefore rank and filter candidates using \emph{touch}. For a table-true inequality $c$ on a finite object set $\mathcal{O}$, define its \emph{sharp set}
\[
  \Sharp(c)\ :=\ \{o\in\mathcal{O}:\ c \text{ holds with equality on } o\},
\]
and its \emph{touch number} as $|\Sharp(c)|$.  Intuitively, a high touch number indicates that the bound is supported by many equality witnesses on the current snapshot.  Nevertheless, even after sorting by touch number, many candidates remain redundant: they introduce no new equality witnesses beyond those already explained by earlier bounds, or they describe the same equality locus in a different syntactic form.  To control this redundancy, we adapt the acceptance logic of Fajtlowicz’s Dalmatian heuristic to the setting where the object collection and candidate list are fixed.  We refer to this one-pass variant as the \emph{static-Dalmatian} heuristic.

The original Dalmatian heuristic is \emph{dynamic}: conjectures are proposed and filtered while the conjecture bank and database evolve over time (for example, as counterexamples are discovered and added). In contrast, static-Dalmatian operates on a single snapshot: a finite family of table-true inequalities generated from a fixed table and then ordered (here, by touch number). Its purpose is not to decide the logical implication among inequalities but to return a compact set of representative bounds, while optionally retaining \emph{aliases}—distinct formulas that identify the same equality locus on the snapshot.

The heuristic maintains a running coverage set $S\subseteq\mathcal{O}$ consisting of objects already touched by accepted representatives. A conjecture $c$ is
accepted as \emph{representative} if it contributes at least one new equality witness,
\[
  \Sharp(c)\setminus S \neq \emptyset,
\]
in which case we update $S\leftarrow S\cup \Sharp(c)$. If, instead, $\Sharp(c)$ coincides with the sharp set of some previously accepted representative, then $c$ is retained as an \emph{alias}, since it provides a different syntactic
description of the same equality structure.

Given a list $\mathcal{C} = (c_1,\dots,c_M)$ sorted in non-increasing order by touch number, the static-Dalmatian filter performs a single scan. It accepts $c_1$ as a representative and initializes $S\leftarrow \Sharp(c_1)$. For each $j=2,\dots,M$, it computes $\Sharp(c_j)$ and: (i) accepts $c_j$ as a representative if $\Sharp(c_j)\setminus S\neq\emptyset$, updating $S\leftarrow S\cup \Sharp(c_j)$;
(ii) otherwise, retains $c_j$ as an alias if $\Sharp(c_j)=\Sharp(c_i)$ for some previously accepted representative $c_i$; and (iii) discards $c_j$ if neither condition holds. The output is the set of accepted conjectures (representatives and aliases).

If $N=|\mathcal{O}|$ and $M$ is the number of candidate conjectures, then a direct implementation that computes $\Sharp(c)$ by scanning all $N$ objects has worst-case runtime $O(MN)$. In practice, constant factors can be reduced by caching computed sharp sets and hashing them to detect aliases efficiently.

\subsection{Full \emph{TxGraffiti} conjecture generation}
\label{subsec:full-generation}

In the previous two subsections, we described the two core components of the baseline pipeline: (i) generating table-true one-predictor affine bounds by solving~\eqref{eq:model1}--\eqref{eq:model2}, and (ii) filtering the resulting candidate list to remove redundancy on the current snapshot via static-Dalmatian. We now summarize how these pieces come together into an end-to-end routine that produces conjectures simultaneously for multiple target invariants.

Fix a snapshot table with numerical invariant columns and Boolean predicate columns.  Let $\mathcal{A}$ be a set of target invariants (columns) for which we wish to conjecture bounds, let $\mathcal{P}$ be the set of available Boolean predicates (hypothesis classes), and let $\mathcal{X}(\alpha)$ denote the admissible numerical predictors for target $\alpha$ (for example, all numerical columns distinct from $\alpha$).  For each target $\alpha\in\mathcal{A}$ and each pair $(x,P)\in\mathcal{X}(\alpha)\times\mathcal{P}$, we solve the upper and lower-bound models~\eqref{eq:model1}--\eqref{eq:model2} restricted to the support set $I(P)$.  Feasible solutions yield affine inequalities
\[
\alpha \le mx+b
\qquad\text{and}\qquad
\alpha \ge mx+b
\]
that are table-true on the predicate-restricted subpopulation.  Each candidate is annotated with its predicate $P$, its touch number, and its sharp set.  Candidates are then sorted by touch number to prioritize bounds that are sharp on many stored instances, and the static-Dalmatian heuristic can be applied to eliminate candidates that add no new sharp witnesses on the current snapshot, while retaining aliases when multiple formulas describe the same equality locus.

Algorithmically, the procedure is a loop over the targets:
\begin{enumerate}
  \item For each target $\alpha\in\mathcal{A}$, enumerate all predictor--predicate pairs $(x,P)\in\mathcal{X}(\alpha)\times\mathcal{P}$ and solve~\eqref{eq:model1}--\eqref{eq:model2} to obtain table-true affine upper and lower bounds.
  \item Compute touch numbers and sharp sets and sort the resulting upper and lower candidate lists by touch number.
  \item Optionally apply static-Dalmatian (Section~\ref{sec:static-dalmatian}) separately to the upper and lower lists.
  \item Merge the surviving conjectures across all targets and sort them for presentation (e.g., by touch number).
\end{enumerate}

This stage provides explicit snapshot-dependent guarantees: every returned inequality is table-true on the selected snapshot table (under its stated predicate), and the optimization objectives in~\eqref{eq:model1}--\eqref{eq:model2} favor bounds with small total slack on the corresponding support set.  Touch number and static-Dalmatian then act as prioritization and de-duplication mechanisms, producing a compact list of representative conjectures together with metadata (sharp instances and aliases) useful for subsequent mathematical scrutiny.  For transparency and reproducibility, the companion repository includes a driver routine that implements this full pipeline and scripts that reproduce the conjecture lists and figures reported in this paper.\footnote{\url{https://github.com/RandyRDavila/Automated-conjecturing-in-mathematics-with-TxGraffiti}}

\section{Example workflow: total domination}
\label{sec:workflow-total-dom}

We illustrate a typical ``\emph{TxGraffiti}-to-theorem'' workflow using the \emph{TxGraffiti} web application.\footnote{\url{https://txgraffiti.streamlit.app}}
The goal is demonstrative rather than exhaustive: we fix one target invariant (the \emph{total domination number} $\gamma_t$), generate candidate bounds under a simple hypothesis, quickly refute several spurious patterns, and then (i) extract one
clean statement that admits a short proof by standard estimates, and (ii) highlight one remaining table-true conjecture that appears to merit further study. We use standard terminology as in West~\cite{West}; see Appendix~\ref{app:terminology} for definitions of $\gamma_t(G)$, $\gamma(G)$, $R(G)$, and graph products $\square$ and $\times$.

\medskip
We begin with a conjecture produced by \emph{TxGraffiti} when restricting the snapshot by the Boolean hypothesis $P(G)\equiv (G\ \text{connected})\wedge(\Delta(G)\le 3)$. Specifically, the system returns the table-true inequality $\gamma_t(G)\ge \tfrac{2}{3}R(G)$ in this restricted snapshot.  Although it first appears as a numerical pattern, it becomes transparent once translated into standard notation: it follows by combining two elementary bounds in terms of $n(G)$ (proved below).

\begin{conj}[\emph{TxGraffiti}]\label{conj:subcubic-gt-Randic}
If $G$ is a connected graph of order $n(G)\ge 2$ with $\Delta(G)\le 3$, then
\[
\gamma_t(G)\ \ge\ \frac{2}{3}\,R(G).
\]
\end{conj}

A common step in our workflow is to test whether a conjecture obtained under a restricted hypothesis is simply a specialization of a more general inequality.  Here, the hypothesis $\Delta(G)\le 3$ suggests replacing the constant $\tfrac{2}{3}$ by the parameter $\tfrac{2}{\Delta(G)}$, and indeed the conjecture above is exactly the specialization of the following $\Delta$-dependent bound at $\Delta(G)=3$.

\begin{thm}\label{thm:gt-Randic-general}
If $G$ is a connected graph of order $n(G)\ge 2$, then
\[
\gamma_t(G)\ \ge\ \frac{2}{\Delta(G)}\,R(G).
\]
\end{thm}

\begin{proof}
Let $D \subseteq V(G)$ be a minimum total dominating set, and so, every vertex in $G$ has a neighbor in $D$. Since every vertex has a
neighbor in $D$, we have $V(G)\subseteq\bigcup_{v\in D}N(v)$ and hence
\[
n(G)\ \le\ \sum_{v\in D}|N(v)|
\ =\ \sum_{v\in D} d(v)
\ \le\ \Delta(G)\,|D|.
\]
Thus, $\gamma_t(G)=|D|\ge n(G)/\Delta(G)$.

\smallskip
\noindent
Next, $G$ has no isolated vertices (because it is connected and $n(G)\ge 2$). For each edge $uv\in E(G)$, AM--GM yields
\[
\frac{1}{\sqrt{d(u)d(v)}}\ \le\ \frac12\!\left(\frac1{d(u)}+\frac1{d(v)}\right).
\]
Summing over edges and collecting by vertex gives
\[
R(G)\ \le\ \frac12\sum_{v\in V(G)} d(v)\cdot\frac1{d(v)}\ =\ \frac{n(G)}{2}.
\]
Combining $R(G)\le n(G)/2$ with $\gamma_t(G)\ge n(G)/\Delta(G)$ yields
\[
\frac{2}{\Delta(G)}\,R(G)\ \le\ \frac{n(G)}{\Delta(G)}\ \le\ \gamma_t(G),
\]
as claimed.
\end{proof}

\begin{corollary}\label{cor:subcubic-gt-Randic}
If $G$ is a connected graph of order $n(G)\ge 2$ and maximum degree $\Delta(G)\le 3$, then
\[
\gamma_t(G)\ \ge\ \frac{2}{3}\,R(G).
\]
In particular, Conjecture~\ref{conj:subcubic-gt-Randic} holds.
\end{corollary}

\begin{proof}
Apply Theorem~\ref{thm:gt-Randic-general} and note that $\Delta(G)\le 3$ implies
$\frac{2}{\Delta(G)}\ge \frac{2}{3}$.
\end{proof}

\medskip
Not all \emph{TxGraffiti} conjectures reduce to such direct combinations of standard estimates.  In our runs, several candidates were quickly refuted by well-known constructions, but a smaller subset remained table-true after basic counterexample searches and did not appear to follow from routine containment or monotonicity arguments. The following conjecture is representative: it compares total domination in the Cartesian product to domination in the direct product.

\begin{conj}[\emph{TxGraffiti}]\label{conj:cart-total-vs-direct-dom}
If $G$ and $H$ are connected graphs with $n(G),n(H)\ge 2$, then
\[
\gamma_t(G \: \square \:H)\ \ge\ \gamma(G\times H).
\]
\end{conj}

Conjecture~\ref{conj:cart-total-vs-direct-dom} is typical of the statements we regard as most promising: it is table-true on the curated datasets used in our runs, it has many sharp instances, and it does not follow from an obvious containment relation between the two products.  Moreover, several superficially similar product inequalities were rapidly refuted, often by exploiting disconnectedness and sparsity phenomena in the direct product, suggesting that any proof here would require genuinely new structural input rather than a routine estimate.

\section{\emph{TxGraffiti}-stimulated results in graph theory}
\label{sec:results}

This section records a small selection of graph-theoretic results whose discovery was \emph{stimulated} by \emph{TxGraffiti}.  By ``stimulated'' we mean that each statement below first appeared as a conjecture output by \emph{TxGraffiti} from a curated snapshot table of graph invariants, and then proceeded through the usual mathematical pipeline—counterexample search, refinement of hypotheses, and proof development—culminating (in the cases collected here) in published theorems or peer-reviewed preprints.

Our aim is not to provide a comprehensive bibliography of \emph{TxGraffiti}-related outcomes. Rather, we document concrete downstream results: specific program-generated conjectures that were subsequently proved in the literature. For ease of reference, each theorem is annotated with the corresponding attribution. Several of these results concern zero forcing and its variants; others relate classical cover- and domination-type invariants. All graph-theoretic terminology and parameters used in the following (including $\alpha,\mu,\beta,\gamma,\gamma_t, Z, Z_t,$ and $Z_{+}$) are defined in the Appendix~\ref{app:terminology}, for standard terminology, we refer the reader to West~\cite{West}.

\subsection{Independence versus matching in regular graphs}
\label{subsec:alpha-vs-mu}

\emph{TxGraffiti} conjectured that in regular graphs, the independence number $\alpha(G)$ is at most the matching number $\mu(G)$.  The parameters $\alpha(G)$ and $\mu(G)$ lie at the core of graph theory: $\alpha(G)$ is a basic extremal measure of sparsity, while $\mu(G)$ is a canonical packing parameter and a cornerstone of structural and algorithmic graph theory.  The resulting clean
inequality
\[
\alpha(G)\le \mu(G)
\]
for $r$-regular graphs ($r>0$) was proved in~\cite{CaDaPe2020}.  From an algorithmic viewpoint, this is also meaningful: computing $\alpha(G)$ is $\NP$-hard in general, whereas $\mu(G)$ is computable in polynomial time (e.g., via Edmonds' matching algorithm). Thus, in regular graphs, Theorem~\ref{thm:txg-alpha-mu-regular} yields a nontrivial efficiently computable upper bound on $\alpha(G)$ in terms of a standard polynomial-time invariant, and it does so with equality on natural sharpness families.

\begin{thm}[\emph{TxGraffiti}; proved in~\cite{CaDaPe2020}]
\label{thm:txg-alpha-mu-regular}
If $G$ is an $r$-regular graph with $r>0$, then
\[
  \alpha(G) \le \mu(G).
\]
\end{thm}

Moreover, \cite{CaDaPe2020} proves a degree-weighted refinement that extends the regular-graph statement to all connected graphs by balancing $\alpha(G)$ and $\mu(G)$ against the minimum and maximum degrees.  The regular case follows immediately by setting $\delta(G)=\Delta(G)=r$.

\begin{thm}[\cite{CaDaPe2020}]
\label{thm:txg-alpha-mu-general}
If $G$ is a connected graph of order $n(G)\ge 2$, then
\[
  \delta(G)\,\alpha(G) \le \Delta(G)\,\mu(G).
\]
\end{thm}

\smallskip
\noindent
Theorem~\ref{thm:txg-alpha-mu-general} interpolates smoothly between highly regular and highly irregular graphs: as the degree spread $\Delta(G)/\delta(G)$ narrows, the matching number increasingly constrains the size of an independent set, recovering Theorem~\ref{thm:txg-alpha-mu-regular} in the regular setting.

\subsection{Zero forcing related results}
\label{subsec:zf-related}

Zero forcing originated as a combinatorial tool for bounding matrix parameters (see~\cite{AIM-Workshop}), but it has also developed into a rich intrinsic graph invariant defined by a local propagation rule.  A consistent theme in the \emph{TxGraffiti} output is that, in structured graph families, the forcing dynamics of $Z(G)$ and its variants are often sharply constrained by classical cover and domination parameters.

\emph{TxGraffiti} conjectured, for instance, that claw-freeness forces the zero forcing number to lie below the vertex cover number.  This inequality relates a static certificate (a minimum vertex cover) to a dynamic propagation process
(zero forcing): in a claw-free graph, the restricted local neighborhood structure allows one to convert an appropriately organized cover into an initial forcing set.  The proof in~\cite{TxGraffiti-2024} is constructive in this sense, producing a zero forcing set from an appropriate vertex cover.

\begin{thm}[\emph{TxGraffiti}; proved in~\cite{TxGraffiti-2024}]
\label{thm:txg-Z-beta-clawfree}
If $G$ is a connected claw-free graph, then
\[
  Z(G) \le \beta(G).
\]
\end{thm}

A second recurring pattern is that within highly structured families—most notably cubic graphs—domination-type parameters exert strong control over zero forcing.  The next two theorems formalize this phenomenon for the standard and total variants.  We emphasize that these results are not routine consequences of
general-purpose techniques: the proofs in~\cite{DaHe21a,DaHe19b} require substantial structural analysis and case work, reflecting the complexity of forcing dynamics even in the restrictive setting of connected cubic graphs.

\begin{thm}[\emph{TxGraffiti}; proved in~\cite{DaHe21a}]
\label{thm:txg-Z-2gamma-cubic}
If $G \neq K_{4}$ is a connected and cubic graph, then
\[
  Z(G) \le 2\,\gamma(G).
\]
\end{thm}

\begin{thm}[\emph{TxGraffiti}; proved in~\cite{DaHe19b}]
\label{thm:txg-Zt-3over2-gammat}
If $G \neq K_{3,3}$ is a connected and cubic graph, then
\[
  Z_t(G) \le \tfrac{3}{2}\,\gamma_t(G).
\]
\end{thm}

Restricting further to claw-free cubic graphs, \emph{TxGraffiti} conjectured a stronger additive domination bound.  The following theorem strengthens what can be guaranteed for the full cubic class by exploiting the additional structure imposed by claw-freeness; the proof in~\cite{DavilaZDom2024} also identifies the extremal graphs (i.e. those achieving equality).

\begin{thm}[\emph{TxGraffiti}; proved in~\cite{DavilaZDom2024}]
\label{thm:txg-Z-gamma-plus-2}
If $G$ is a connected, cubic, and claw-free graph, then
\[
  Z(G) \le \gamma(G) + 2,
\]
and this bound is sharp if and only if $G$ is a diamond-necklace or $G$ is the complete graph $K_4$.
\end{thm}

The preceding theorems motivate a broader conjectural relationship between the independence number and the zero forcing number in graphs of maximum degree at most three. This conjecture has been studied extensively and remains a central open problem originating from \emph{TxGraffiti} output.

\begin{conj}[\emph{TxGraffiti}]
\label{conj:txgraffiti-main}
If $G \neq K_4$ is a connected graph with $\Delta(G) \le 3$, then
\[
  Z(G) \le \alpha(G) + 1.
\]
\end{conj}

Conjecture~\ref{conj:txgraffiti-main} has been proved in important special cases, including cubic claw-free graphs~\cite{DaHe19c}.  More recently, Schuerger, Warnberg, and Young show that \emph{almost all} cubic graphs satisfy the weaker bound $Z(G)\le \alpha(G)+2$, and they exhibit an infinite family of cubic graphs
with $Z(G)=\alpha(G)+1$~\cite{SchuergerWarnbergYoung2024Zalpha}.

\smallskip
\noindent
\emph{TxGraffiti} also suggested that claw-free structure collapses the gap between standard and positive semidefinite zero forcing.  Since $Z_{+}(G)$ is motivated by minimum-rank problems for positive semidefinite matrices, the equality below provides a bridge between a forbidden-induced-subgraph condition and an algebraically motivated forcing parameter. The proof in \cite{DavilaSchuergerSmall2024ClawFreeZeqZplus} establishes equality for all claw-free graphs and yields a characterization in terms of induced subgraphs.

\begin{thm}[\emph{TxGraffiti}; proved in~\cite{DavilaSchuergerSmall2024ClawFreeZeqZplus}]
\label{thm:txg-Z-eq-Zplus}
If $G$ is claw-free, then
\[
  Z(G)=Z_{+}(G).
\]
\end{thm}

\begin{corollary}[\cite{DavilaSchuergerSmall2024ClawFreeZeqZplus}]
\label{thm:txg-Z-eq-Zplus}
A connected graph $G$ is claw-free if and only if $Z(H) = Z_{+}(H)$ for every induced subgraph $H \subseteq G$. 
\end{corollary}

\smallskip
\noindent
Taken together, these results illustrate a typical ``\emph{TxGraffiti}-to-theorem'' trajectory: empirically discovered inequalities are subjected to systematic counterexample search, refined to the appropriate hypothesis class (often expressed via forbidden induced subgraphs or degree constraints), and then proved—frequently alongside sharpness analyses, identification of extremal families, and isolation of a small number of exceptional graphs.

\section{Conclusion}
\label{sec:conclusion}

\emph{TxGraffiti}, as described in this manuscript, frames automated conjecturing around a finite versioned \emph{snapshot table} of objects equipped with precomputed numerical invariants and Boolean predicates.  In this setting, conjecturing becomes a sequential optimization-and-filtering pipeline: for a chosen target invariant and hypothesis predicate, the system searches over (predicate, predictor) combinations, solves a sequence of small linear programs to fit coefficients within simple templates (e.g., univariate affine bounds) and outputs \emph{table-true} conditional inequalities of the form
\[
P(G)\ \Rightarrow\ \alpha(G)\le f(G)
\qquad\text{or}\qquad
P(G)\ \Rightarrow\ \alpha(G)\ge f(G),
\]
which can be interpreted as conjectures and then tested, refined, or pursued for proof beyond the snapshot.  Heuristics such as touch/sharp statistics and redundancy filtering are subsequently applied to compress the resulting candidate set into a small collection of informative and human-readable statements.

The account presented here reflects the system in its formative stage and emphasizes the design choice to prioritize fast, interpretable conjectures from fixed evidence tables. Ongoing development has extended the platform’s capabilities, and a modern open-source implementation is available as the \texttt{txgraffiti} Python package on PyPI~\cite{txgraffiti_pypi}.  Future directions include tighter integration between interactive table construction and conjecturing, richer but still interpretable template families (including schema-level feature engineering and multivariate bounds), and collaborative workflows for curating snapshot tables and theory libraries.

\section*{Acknowledgements}

The author gratefully acknowledges the referees for their careful reading and many insightful suggestions, which substantially improved the quality and clarity of this manuscript.

\section*{Data availability}
All data generated or analyzed during this study are included in this published article and its supplementary information files.

\section*{Code availability}
All code used in this paper can be found at this GitHub repository: \url{https://github.com/RandyRDavila/Automated-conjecturing-in-mathematics-with-TxGraffiti}

\section*{Declarations}
The authors declare no competing interests.

\appendix
\section{Graph-theoretic notation and terminology}
\label{app:terminology}

This appendix provides notation and graph-theoretic parameters used in this paper. Unless stated otherwise, all graphs are finite, simple, and
undirected. For a graph $G$, we write $V(G)$ and $E(G)$ for its vertex and edge sets, and we set $n(G)=|V(G)|$ and $m(G)=|E(G)|$.  For $v\in V(G)$, the (open) neighborhood is $N(v)=\{u\in V(G): uv\in E(G)\}$, the degree is $d(v)=|N(v)|$, and we write $\Delta(G)=\max_{v\in V(G)} d(v)$ and $\delta(G)=\min_{v\in V(G)} d(v)$ for the maximum and minimum degrees.

Several classical invariants appear throughout.  An \emph{independent set} is a set $S\subseteq V(G)$ in which no two vertices are adjacent, and the
\emph{independence number} $\alpha(G)$ is the maximum size of such a set.  A \emph{matching} is a set of pairwise edge-disjoint edges, and the
\emph{matching number} $\mu(G)$ is the maximum size of a matching.  A
\emph{vertex cover} is a set $C\subseteq V(G)$ meeting every edge (equivalently, for every $uv\in E(G)$ at least one of $u$ or $v$ lies in $C$), and the \emph{vertex cover number} $\beta(G)$ is the minimum size of a vertex cover. A \emph{dominating set} is a set $D\subseteq V(G)$ such that every vertex in $V(G)\setminus D$ has a neighbor in $D$, and the \emph{domination number} $\gamma(G)$ is the minimum size of a dominating set.  Total domination is the variant in which all vertices must be dominated by a \emph{neighbor}: a graph admits a \emph{total dominating set} precisely when it has no isolated vertices, and in that case a set $D\subseteq V(G)$ is a \emph{total dominating set} if every vertex has a neighbor in $D$.  The \emph{total domination number} $\gamma_t(G)$ is the minimum size of a total dominating set.

We also refer to certain graph classes.  A graph is \emph{$r$-regular} if $d(v)=r$ for all $v\in V(G)$; in particular, it is \emph{cubic} if it is
$3$-regular.  A graph is \emph{claw-free} if it contains no induced subgraph isomorphic to $K_{1,3}$. The \emph{Randi\'c index} of $G$ is
\[
R(G)\;=\;\sum_{uv\in E(G)}\frac{1}{\sqrt{d(u)d(v)}}.
\]

Two standard graph products are used.  The \emph{Cartesian product} $G \: \square \: H$ has vertex set $V(G)\times V(H)$, where $(g,h)$ is adjacent to $(g',h')$ if either $g=g'$ and $hh'\in E(H)$, or $h=h'$ and $gg'\in E(G)$.  The \emph{direct product} (also called the \emph{tensor product}) $G\times H$ likewise has vertex set $V(G)\times V(H)$, but now $(g,h)$ is adjacent to $(g',h')$ precisely when $gg'\in E(G)$ and $hh'\in E(H)$.

Finally, we use zero forcing and two common variants.  In the standard
\emph{zero forcing process}~\cite{AIM-Workshop}, one starts from an initial set $S\subseteq V(G)$ of blue vertices (all others white) and repeatedly applies the color-change rule: a blue vertex $u$ with exactly one white neighbor $w$ forces that neighbor to become blue (written $u\to w$).  A set $S$ is a \emph{zero forcing set} if this process eventually colors all vertices blue, and the \emph{zero forcing number} $Z(G)$ is the minimum size of such a set.  The \emph{total zero forcing number} $Z_t(G)$ is defined analogously, with the additional requirement that the initial blue set induces no isolated vertices (equivalently, every initially blue vertex has a blue neighbor), see~\cite{DavilaThesis2019}.

The \emph{positive semidefinite zero forcing number} $Z_{+}(G)$ is defined using the positive semidefinite color-change rule~\cite{Barioli2010}.  Given a current set $B$ of blue vertices, let $W_1,\dots,W_k$ be the vertex sets of the connected components of the induced subgraph $G[V(G)\setminus B]$.  If $u\in B$ and $w\in W_i$ is the \emph{only} white neighbor of $u$ in the induced subgraph $G[W_i\cup B]$, then $u$ forces $w$ to become blue.  A set $S$ is a \emph{positive semidefinite zero forcing set} if repeated application of this rule colors all vertices blue, and $Z_{+}(G)$ is the minimum size of such a set.

\bibliography{sn-bibliography}

@incollection{Graffitipc,
  author    = {DeLaVi{\~n}a, Ermelinda},
  title     = {Graffiti.pc: {A} variant of {G}raffiti},
  booktitle = {Graphs and Discovery},
  editor    = {Fajtlowicz, Siemion and Fowler, Patrick W. and Hansen, Pierre and Janowitz, Melvin F. and Roberts, Fred S.},
  series    = {DIMACS Series in Discrete Mathematics and Theoretical Computer Science},
  volume    = {69},
  publisher = {American Mathematical Society},
  address   = {Providence, RI},
  year      = {2005},
  pages     = {71--80},
  doi       = {10.1090/dimacs/069}
}

@incollection{GraffitiD,
  author    = {DeLaVi{\~n}a, Ermelinda},
  title     = {Some history of the development of {G}raffiti},
  booktitle = {Graphs and Discovery},
  editor    = {Fajtlowicz, Siemion and Fowler, Patrick W. and Hansen, Pierre and Janowitz, Melvin F. and Roberts, Fred S.},
  series    = {DIMACS Series in Discrete Mathematics and Theoretical Computer Science},
  volume    = {69},
  publisher = {American Mathematical Society},
  address   = {Providence, RI},
  year      = {2005},
  pages     = {81--118},
  doi       = {10.1090/dimacs/069}
}

@incollection{Fajtlowicz-Buckminsterfullerene-2005,
  author    = {Fajtlowicz, Siemion},
  title     = {On representation and characterization of buckminsterfullerene {$C_{60}$}},
  booktitle = {Graphs and Discovery},
  editor    = {Fajtlowicz, Siemion and Fowler, Patrick W. and Hansen, Pierre and Janowitz, Melvin F. and Roberts, Fred S.},
  series    = {DIMACS Series in Discrete Mathematics and Theoretical Computer Science},
  volume    = {69},
  publisher = {American Mathematical Society},
  address   = {Providence, RI},
  year      = {2005},
  pages     = {127--135},
  doi       = {10.1090/dimacs/069}
}

@article{Larson,
  author  = {Larson, Craig E. and Van Cleemput, Nicolas},
  title   = {Automated conjecturing {I}: {F}ajtlowicz's {D}almatian heuristic revisited},
  journal = {Artificial Intelligence},
  volume  = {231},
  pages   = {17--38},
  year    = {2016},
  month   = {2},
  doi     = {10.1016/j.artint.2015.10.002}
}

@article{LarsonVC2017,
  author  = {Larson, C. E. and Van Cleemput, N.},
  title   = {Automated conjecturing {III}: Property-relations conjectures},
  journal = {Annals of Mathematics and Artificial Intelligence},
  volume  = {81},
  pages   = {315--327},
  year    = {2017},
  doi     = {10.1007/s10472-017-9559-5}
}

@article{Robbins,
  author  = {McCune, William},
  title   = {Solution of the Robbins Problem},
  journal = {Journal of Automated Reasoning},
  volume  = {19},
  number  = {3},
  pages   = {263--276},
  year    = {1997},
  doi     = {10.1023/A:1005843212881}
}

@article{SimonNewell,
  author  = {Simon, Herbert A. and Newell, Allen},
  title   = {Heuristic Problem Solving: The Next Advance in Operations Research},
  journal = {Operations Research},
  volume  = {6},
  number  = {1},
  pages   = {1--10},
  year    = {1958},
  doi     = {10.1287/opre.6.1.1}
}

@incollection{Turing,
  author    = {Turing, Alan M.},
  title     = {Intelligent Machinery},
  booktitle = {The Essential Turing: Seminal Writings in Computing, Logic, Philosophy, Artificial Intelligence, and Artificial Life: Plus The Secrets of Enigma},
  editor    = {Copeland, B. Jack},
  publisher = {Oxford University Press},
  year      = {2004},
  pages     = {395--432},
  address = {Oxford},
  note      = {Written in 1948; reprinted in this volume.}
}

@article{WangProgramII,
  author  = {Wang, Hao},
  title   = {Proving Theorems by Pattern Recognition I},
  journal = {Communications of the ACM},
  volume  = {3},
  number  = {4},
  pages   = {220--234},
  year    = {1960},
  doi     = {10.1145/367177.367224}
}

@article{Fajtlowicz-DM-1988,
  author  = {Fajtlowicz, Siemion},
  title   = {On Conjectures of {G}raffiti},
  journal = {Discrete Math.},
  volume  = {72},
  number  = {1-3},
  pages   = {113--118},
  year    = {1988},
  doi     = {10.1016/0012-365X(88)90199-9}
}

@article{Fajtlowicz-1987,
  author  = {Fajtlowicz, Siemion},
  title   = {On Conjectures of {G}raffiti, {II}},
  journal = {Congressus Numerantium},
  volume  = {60},
  pages   = {187--197},
  year    = {1987}
}

@article{Fajtlowicz-III-1988,
  author  = {Fajtlowicz, Siemion},
  title   = {On Conjectures of {G}raffiti, {III}},
  journal = {Congressus Numerantium},
  volume  = {66},
  pages   = {23--32},
  year    = {1988}
}

@article{Fajtlowicz-IV-1990,
  author  = {Fajtlowicz, Siemion},
  title   = {On Conjectures of {G}raffiti, {IV}},
  journal = {Congressus Numerantium},
  volume  = {70},
  pages   = {231--240},
  year    = {1990}
}

@inproceedings{Fajtlowicz-V-1995,
  author    = {Fajtlowicz, Siemion},
  title     = {On Conjectures of {G}raffiti, {V}},
  booktitle = {Graph Theory, Combinatorics, and Algorithms: Proceedings of the Seventh Quadrennial International Conference on the Theory and Applications of Graphs},
  editor    = {Alavi, Yousef and Schwenk, Allen},
  volume    = {1},
  pages     = {367--376},
  publisher = {John Wiley \& Sons},
  address   = {New York},
  year      = {1995}
}

@inproceedings{Fajtlowicz-Clemson-1989,
  author       = {Fajtlowicz, Siemion},
  title        = {On Conjectures and Methods of {G}raffiti},
  booktitle    = {Proceedings of the Fourth Clemson Mini-Conference on Graph Theory and Combinatorics},
  address      = {Clemson, SC},
  year         = {1989},
  note         = {Clemson University}
}

@article{Fajtlowicz-Larson-2003,
  author  = {Fajtlowicz, Siemion and Larson, Craig E.},
  title   = {Graph-theoretic independence as a predictor of fullerene stability},
  journal = {Chemical Physics Letters},
  volume  = {377},
  pages   = {485--490},
  year    = {2003},
  doi     = {10.1016/S0009-2614(03)01133-3}
}

@article{Fajtlowicz-John-Sachs-2005,
  author  = {Fajtlowicz, Siemion and John, Peter E. and Sachs, Horst},
  title   = {On maximum matchings and eigenvalues of benzenoid graphs},
  journal = {Croatica Chemica Acta},
  volume  = {78},
  number  = {2},
  pages   = {195--201},
  year    = {2005},
  url     = {https://hrcak.srce.hr/en/12}
}

@article{Fajtlowicz-Waller-1986,
  author  = {Fajtlowicz, Siemion and Waller, William A.},
  title   = {On Two Conjectures of {G}raffiti},
  journal = {Congressus Numerantium},
  volume  = {55},
  pages   = {51--56},
  year    = {1986}
}

@article{Chung-1988,
  author  = {Chung, F. R. K.},
  title   = {The Average Distance and the Independence Number},
  journal = {Journal of Graph Theory},
  volume  = {12},
  pages   = {229--235},
  year    = {1988},
  doi     = {10.1002/jgt.3190120213}
}

@article{Alon-Seymour-1989,
  author  = {Alon, Noga and Seymour, Paul D.},
  title   = {A Counterexample to the Rank-Coloring Conjecture},
  journal = {Journal of Graph Theory},
  volume  = {13},
  number  = {4},
  pages   = {523--525},
  year    = {1989},
  doi     = {10.1002/jgt.3190130413}
}

@techreport{Beezer-1989,
  author      = {Beezer, Robert A. and Riegsecker, J. and Smith, B. A.},
  title       = {On Conjectures of {G}raffiti Concerning Regular Graphs},
  institution = {University of Puget Sound, Department of Mathematics and Computer Science},
  number      = {89-1},
  year        = {1989},
  month       = {8}
}

@article{Favaron-1990,
  author  = {Favaron, Odile and Mah{\'e}o, Michel and Sacl{\'e}, J.-F.},
  title   = {Some results on conjectures of {G}raffiti {I}},
  journal = {Ars Combinatoria},
  volume  = {29C},
  pages   = {90--106},
  year    = {1990}
}

@techreport{Favaron-Research-1991,
  author      = {Favaron, Odile and Mah{\'e}o, Michel and Sacl{\'e}, J.-F.},
  title       = {Some results on conjectures of {G}raffiti {III}},
  institution = {LRI, Universit{\'e} de Paris-Sud},
  number      = {670},
  type        = {Research Report},
  year        = {1991}
}

@article{Favaron-1991,
  author  = {Favaron, Odile and Mah{\'e}o, Maryvonne and Sacl{\'e}, Jean-Fran{\c{c}}ois},
  title   = {On the residue of a graph},
  journal = {Journal of Graph Theory},
  volume  = {15},
  number  = {1},
  pages   = {39--64},
  year    = {1991},
  doi     = {10.1002/jgt.3190150107}
}

@article{Fajtlowicz-McColgan-1995,
  author  = {Fajtlowicz, Siemion and McColgan, Tamara and Reid, Talmage and Staton, William},
  title   = {Ramsey Numbers for Induced Regular Subgraphs},
  journal = {Ars Combinatoria},
  volume  = {39},
  pages   = {149--154},
  year    = {1995},
  url     = {https://combinatorialpress.com/ars-articles/volume-039-ars-articles/ramsey-numbers-for-induced-regular-subgraphs/}
}

@article{Favaron-2003,
  author  = {Favaron, Odile and Mah{\'e}o, Maryvonne and Sacl{\'e}, Jean-Fran{\c{c}}ois},
  title   = {The Randi{\'c} index and other {G}raffiti parameters of graphs},
  journal = {MATCH Commun. Math. Comput. Chem.},
  volume  = {47},
  pages   = {7--23},
  year    = {2003},
  url     = {https://match.pmf.kg.ac.rs/content47.htm}
}

@article{Favaron-DM-1993,
  author  = {Favaron, Odile and Mah{\'e}o, Maryvonne and Sacl{\'e}, Jean-Fran{\c{c}}ois},
  title   = {Some eigenvalue properties in graphs (conjectures of {G}raffiti {II})},
  journal = {Discrete Math.},
  volume  = {111},
  number  = {1-3},
  pages   = {197--220},
  year    = {1993},
  doi     = {10.1016/0012-365X(93)90156-N}
}

@article{Beezer-Riegsecker-Smith-2001,
  author  = {Beezer, Robert A. and Riegsecker, John E. and Smith, Bryan A.},
  title   = {Using minimum degree to bound average distance},
  journal = {Discrete Math.},
  volume  = {226},
  number  = {1-3},
  pages   = {365--371},
  year    = {2001},
  doi     = {10.1016/S0012-365X(00)00156-4},
  url     = {https://doi.org/10.1016/S0012-365X(00)00156-4}
}

@article{Caro-1998,
  author  = {Caro, Yair},
  title   = {Colorability, frequency and {Graffiti--119}},
  journal = {J. Comb. Math. Comb. Comput.},
  volume  = {27},
  pages   = {129--134},
  year    = {1998}
}

@article{Codenotti-2000,
  author  = {Codenotti, Bruno and Del Corso, Gianna and Manzini, Giovanni},
  title   = {Matrix rank and communication complexity},
  journal = {Linear Algebra Appl.},
  volume  = {304},
  number  = {1-3},
  pages   = {193--200},
  year    = {2000},
  doi     = {10.1016/S0024-3795(99)00226-8},
  url     = {https://doi.org/10.1016/S0024-3795(99)00226-8}
}

@article{Cygan-2012,
  author  = {Cygan, Marek and Pilipczuk, Micha{\l} and {\v{S}}krekovski, Riste},
  title   = {On the inequality between radius and Randi{\'c} index for graphs},
  journal = {MATCH Commun. Math. Comput. Chem.},
  volume  = {67},
  number  = {2},
  pages   = {451--466},
  year    = {2012},
  url     = {https://match.pmf.kg.ac.rs/electronic_versions/Match67/n2/match67n2_451-466.pdf}
}

@article{Lenat_1,
  author  = {Lenat, Douglas B.},
  title   = {The ubiquity of discovery},
  journal = {Artificial Intelligence},
  volume  = {9},
  number  = {3},
  pages   = {257--285},
  year    = {1977},
  doi     = {10.1016/0004-3702(77)90024-8}
}

@incollection{Lenat_2,
  author    = {Lenat, Douglas B.},
  title     = {On automated scientific theory formation: {A} case study using the {AM} program},
  booktitle = {Machine Intelligence 9},
  editor    = {Hayes, John E. and Michie, Donald and Mikulich, L. I.},
  publisher = {Halsted Press},
  address   = {New York},
  pages     = {251--283},
  year      = {1979}
}

@article{Lenat_3,
  author  = {Lenat, Douglas B.},
  title   = {The nature of heuristics},
  journal = {Artificial Intelligence},
  volume  = {19},
  number  = {2},
  pages   = {189--249},
  year    = {1982},
  month   = {10},
  doi     = {10.1016/0004-3702(82)90036-4}
}

@article{graphedron_1,
  author  = {M{\'e}lot, Hadrien},
  title   = {Facet defining inequalities among graph invariants: The system {GraPHedron}},
  journal = {Discrete Applied Mathematics},
  volume  = {156},
  number  = {10},
  pages   = {1875--1891},
  year    = {2008},
  doi     = {10.1016/j.dam.2007.09.005}
}

@article{Dankelmann-Dlamini-Swart-2005,
  author  = {Dankelmann, Peter and Dlamini, Gcina and Swart, Henda C.},
  title   = {Upper bounds on distance measures in {$K_{3,3}$}-free graphs},
  journal = {Utilitas Mathematica},
  volume  = {67},
  pages   = {205--222},
  year    = {2005}
}

@article{Dankelmann-Swart-Oellermann-1998,
  author  = {Dankelmann, Peter and Swart, Henda C. and Oellermann, Ortrud R.},
  title   = {On three conjectures of {G}raffiti},
  journal = {Journal of Combinatorial Mathematics and Combinatorial Computing},
  volume  = {26},
  pages   = {131--137},
  year    = {1998},
  url     = {https://combinatorialpress.com/jcmcc-articles/volume-026/on-three-conjectures-of-graffiti/}
}

@article{Bollobas-Erdos-1998,
  author  = {Bollob{\'a}s, B{\'e}la and Erd{\H o}s, Paul},
  title   = {Graphs of Extremal Weights},
  journal = {Ars Combinatoria},
  volume  = {50},
  pages   = {225--233},
  year    = {1998},
  url     = {https://combinatorialpress.com/ars-articles/volume-050-ars-articles/graphs-of-extremal-weights/}
}

@article{Bollobas-Riordan-1998,
  author  = {Bollob{\'a}s, B{\'e}la and Riordan, Oliver M.},
  title   = {On Some Conjectures of {G}raffiti},
  journal = {Discrete Mathematics},
  volume  = {179},
  number  = {1-3},
  pages   = {223--230},
  year    = {1998},
  doi     = {10.1016/S0012-365X(97)00093-9}
}

@article{Epstein_2,
  author  = {Epstein, Susan L.},
  title   = {Learning and discovery: One system's search for mathematical knowledge},
  journal = {Computational Intelligence},
  volume  = {4},
  number  = {1},
  pages   = {42--53},
  year    = {1988},
  url     = {https://www.cs.hunter.cuny.edu/~epstein/papers/Epstein-1988-Computational_Intelligence.pdf}
}

@incollection{devillez2019,
  author    = {Devillez, Gauvain and Hauweele, Pierre and M{\'e}lot, Hadrien},
  title     = {{PHOEG} Helps to Obtain Extremal Graphs},
  booktitle = {Operations Research Proceedings 2018: Selected Papers of the Annual International Conference of the German Operations Research Society (GOR), Brussels, Belgium, September 12--14, 2018},
  editor    = {Fortz, Bernard and Labb{\'e}, Martine},
  series    = {Operations Research Proceedings},
  publisher = {Springer},
  address   = {Cham},
  year      = {2019},
  pages     = {251--257},
  doi       = {10.1007/978-3-030-18500-8_32}
}

@article{AGX_1,
  author  = {Caporossi, G. and Hansen, P.},
  title   = {Variable neighborhood search for extremal graphs: {I}. The {AutoGraphiX} system},
  journal = {Discrete Mathematics},
  volume  = {212},
  number  = {1--2},
  pages   = {29--44},
  year    = {2000},
  doi     = {10.1016/S0012-365X(99)00206-X},
}

@article{AGX_2,
  author  = {Caporossi, G. and Hansen, P.},
  title   = {Variable neighborhood search for extremal graphs: {V}. Three ways to automate finding conjectures},
  journal = {Discrete Mathematics},
  volume  = {276},
  number  = {1--3},
  pages   = {81--94},
  year    = {2004},
  doi     = {10.1016/S0012-365X(03)00311-X},
}

@article{Hansen-2009,
  author  = {Hansen, P. and Hertz, A. and Kilani, R. and Marcotte, O. and Schindl, D.},
  title   = {Average distance and maximum induced forest},
  journal = {J. Graph Theory},
  volume  = {60},
  number  = {1},
  pages   = {31--54},
  year    = {2009},
  doi     = {10.1002/jgt.20344},
}

@article{Jelen-1999,
  author  = {Jelen, F.},
  title   = {$k$-Independence and the $k$-residue of a graph},
  journal = {J. Graph Theory},
  volume  = {32},
  number  = {3},
  pages   = {241--249},
  year    = {1999},
  doi     = {10.1002/(SICI)1097-0118(199911)32:3<241::AID-JGT4>3.0.CO;2-S},
}

@article{Firby-1997,
  author  = {Firby, Peter and Haviland, Julie},
  title   = {Independence and average distance in graphs},
  journal = {Discrete Applied Mathematics},
  volume  = {75},
  number  = {1},
  pages   = {27--37},
  year    = {1997},
  doi     = {10.1016/S0166-218X(96)00078-9}
}

@article{Griggs-Kleitman-1994,
  author  = {Griggs, Jerrold R. and Kleitman, Daniel J.},
  title   = {Independence and the {H}avel--{H}akimi residue},
  journal = {Discrete Mathematics},
  volume  = {127},
  number  = {1--3},
  pages   = {209--212},
  year    = {1994},
  doi     = {10.1016/0012-365X(92)00479-B}
}

@article{Wang-1997,
  author  = {Wang, Liuxing},
  title   = {On one of {G}raffiti's conjecture (583)},
  journal = {Applied Mathematics and Mechanics},
  volume  = {18},
  pages   = {381--383},
  year    = {1997},
  month   = {4},
  doi     = {10.1007/BF02457552}
}

@article{Zhang-2004,
  author  = {Zhang, X.-D.},
  title   = {On the two conjectures of {G}raffiti},
  journal = {Linear Algebra and its Applications},
  volume  = {385},
  pages   = {369--379},
  year    = {2004},
  doi     = {10.1016/j.laa.2003.12.014}
}

@article{Doslic-Reti-2011,
  author  = {Do{\v{s}}li{\'c}, Tomislav and Reti{\'c}, Tam{\'a}s},
  title   = {Spectral properties of fullerene graphs},
  journal = {MATCH Communications in Mathematical and in Computer Chemistry},
  volume  = {66},
  pages   = {733--742},
  year    = {2011},
  url     = {https://match.pmf.kg.ac.rs/electronic_versions/Match66/n3/match66_733-742.pdf}
}

@article{Fowler-1997,
  author  = {Fowler, P. W.},
  title   = {Fullerene graphs with more negative than positive eigenvalues: The exceptions that prove the rule of electron deficiency?},
  journal = {J. Chem. Soc., Faraday Trans.},
  volume  = {93},
  pages   = {1--3},
  year    = {1997},
  doi     = {10.1039/A605413G}
}

@article{Fowler-1998,
  author  = {Fowler, P. W. and Hansen, P. and Rogers, K. M. and Fajtlowicz, S.},
  title   = {{$C_{60}Br_{24}$} as a chemical illustration of graph theoretical independence},
  journal = {J. Chem. Soc., Perkin Trans. 2},
  pages   = {1531--1534},
  year    = {1998},
  doi     = {10.1039/A803459A}
}

@incollection{Fowler-1999,
  author    = {Fowler, Patrick W. and Rogers, Kevin M. and Fajtlowicz, Siemion and Hansen, Pierre and Caporossi, Gilles},
  title     = {Facts and Conjectures about Fullerene Graphs: Leapfrog, Cylinder and Ramanujan Fullerenes},
  booktitle = {Algebraic Combinatorics and Applications},
  booksubtitle = {Proceedings of the Euroconference, Algebraic Combinatorics and Applications (ALCOMA), held in G{\"o}{\ss}weinstein, Germany, September 12--19, 1999},
  editor    = {Betten, Anton and Kohnert, Axel and Laue, Reinhard and Wassermann, Alfred},
  publisher = {Springer},
  address   = {Berlin, Heidelberg},
  year      = {2001},
  pages     = {134--146},
  doi       = {10.1007/978-3-642-59448-9_10}
}

@article{Raayoni2021RamanujanMachine,
  title   = {Generating conjectures on fundamental constants with the {Ramanujan} Machine},
  author  = {Raayoni, Gal and Gottlieb, Shahar and Burshtein, Noam and Kaminer, Ido and others},
  journal = {Nature},
  volume  = {590},
  number  = {7844},
  pages   = {67--73},
  year    = {2021},
  doi     = {10.1038/s41586-021-03229-4}
}

@article{CaDaPe2020,
  author  = {Caro, Yair and Davila, Randy and Pepper, Ryan},
  title   = {New results relating matching and independence},
  journal = {Discussiones Mathematicae Graph Theory},
  volume  = {42},
  number  = {3},
  pages   = {921--935},
  year    = {2020},
  doi     = {10.7151/dmgt.2317}
}

@misc{DavilaZDom2024,
  author       = {Davila, Randy},
  title        = {Another conjecture of {TxGraffiti} concerning zero forcing and domination in graphs},
  year         = {2024},
  eprint       = {2406.19231},
  eprinttype   = {arXiv},
  primaryclass = {math.CO},
  doi          = {10.48550/arXiv.2406.19231},
  url          = {https://arxiv.org/abs/2406.19231},
  note         = {Under revision (Electronic Journal of Combinatorics)}
}

@phdthesis{DavilaThesis2019,
  author = {Davila, Randy},
  title  = {Total and Zero Forcing in Graphs},
  school = {University of Johannesburg},
  year   = {2019},
  type   = {D.Phil. thesis},
  url    = {https://hdl.handle.net/10210/295628}
}

@article{DaHe21a,
  author  = {Davila, Randy and Henning, Michael A.},
  title   = {Zero forcing versus domination in cubic graphs},
  journal = {Journal of Combinatorial Optimization},
  volume  = {41},
  pages   = {553--577},
  year    = {2021},
  doi     = {10.1007/s10878-020-00692-z}
}

@article{DaHe19b,
  author  = {Davila, Randy and Henning, Michael A.},
  title   = {Total forcing versus total domination in cubic graphs},
  journal = {Applied Mathematics and Computation},
  volume  = {354},
  pages   = {385--395},
  year    = {2019},
  doi     = {10.1016/j.amc.2019.02.060}
}

@article{DaHe19c,
  author  = {Davila, Randy and Henning, Michael A.},
  title   = {Zero forcing in claw-free cubic graphs},
  journal = {Bulletin of the Malaysian Mathematical Sciences Society},
  volume  = {43},
  pages   = {673--688},
  year    = {2020},
  doi     = {10.1007/s40840-018-00705-5}
}

@article{TxGraffiti-2024,
  author  = {Brimkov, Boris and Davila, Randy and Schuerger, Houston and Young, Michael},
  title   = {On a conjecture of {TxGraffiti}: Relating zero forcing and vertex covers in graphs},
  journal = {Discrete Applied Mathematics},
  volume  = {359},
  pages   = {290--302},
  year    = {2024},
  doi     = {10.1016/j.dam.2024.08.006}
}

@inproceedings{Colton_1,
  author    = {Colton, Simon and Bundy, Alan and Walsh, Toby},
  title     = {Automated concept formation in pure mathematics},
  booktitle = {Proceedings of the Sixteenth International Joint Conference on Artificial Intelligence (IJCAI-99), Stockholm, Sweden, July 31--August 6, 1999},
  editor    = {Dean, Thomas},
  volume    = {2},
  pages     = {786--791},
  publisher = {Morgan Kaufmann},
  address   = {San Francisco, CA},
  year      = {1999},
  url       = {https://www.ijcai.org/Proceedings/99-2/Papers/018.pdf}
}

@article{Colton_2,
  author  = {Colton, Simon},
  title   = {Refactorable numbers---a machine invention},
  journal = {J. Integer Seq.},
  volume  = {2},
  pages   = {Article 99.1.2},
  year    = {1999},
  url     = {https://cs.uwaterloo.ca/journals/JIS/colton/joisol.html}
}

@book{Colton_3,
  author    = {Colton, Simon},
  title     = {Automated Theory Formation in Pure Mathematics},
  series    = {Distinguished Dissertations},
  publisher = {Springer London},
  address   = {London},
  year      = {2002},
  doi       = {10.1007/978-1-4471-0147-5},
  isbn      = {978-1-85233-609-7}
}

@article{AIM-Workshop,
  author  = {{AIM Minimum Rank--Special Graphs Work Group}},
  title   = {Zero forcing sets and the minimum rank of graphs},
  journal = {Linear Algebra and its Applications},
  volume  = {428},
  number  = {7},
  pages   = {1628--1648},
  year    = {2008},
  doi     = {10.1016/j.laa.2007.10.009}
}

@book{West,
  author    = {West, Douglas B.},
  title     = {Introduction to Graph Theory},
  edition   = {2},
  publisher = {Prentice Hall},
  address   = {Upper Saddle River, NJ},
  year      = {2001},
  isbn      = {0-13-014400-2}
}

@article{Yue-2020,
  author  = {Yue, Jun and Zhang, Shizhen and Zhu, Yiping and Klav{\v{z}}ar, Sandi and Shi, Yongtang},
  title   = {The annihilation number does not bound the $2$-domination number from the above},
  journal = {Discrete Mathematics},
  volume  = {343},
  number  = {6},
  pages   = {111707},
  year    = {2020},
  doi     = {10.1016/j.disc.2019.111707}
}

@article{Davies2021AdvancingMathAI,
  author  = {Davies, Alex and Veli{\v{c}}kovi{\'c}, Petar and Buesing, Lars and Blackwell, Sam and Zheng, Daniel and Toma{\v{s}}ev, Nenad and Tanburn, Richard and Battaglia, Peter W. and Blundell, Charles and Juh{\'a}sz, Andr{\'a}s and Lackenby, Marc and Williamson, Geordie and Hassabis, Demis and Kohli, Pushmeet},
  title   = {Advancing mathematics by guiding human intuition with {AI}},
  journal = {Nature},
  volume  = {600},
  number  = {7887},
  pages   = {70--74},
  year    = {2021},
  doi     = {10.1038/s41586-021-04086-x},
}

@article{DoughertyBlissZeilberger2023,
  author  = {Dougherty-Bliss, Robert and Zeilberger, Doron},
  title   = {Automatic conjecturing and proving of exact values of some infinite families of infinite continued fractions},
  journal = {The Ramanujan Journal},
  volume  = {61},
  pages   = {31--47},
  year    = {2023},
  doi     = {10.1007/s11139-020-00345-z},
  note    = {Published online 24 Feb 2021},
}

@misc{Yamamoto2024RamanujanMachineProof,
  author       = {Yamamoto, Shuma},
  title        = {Proof and generalization of conjectures of {R}amanujan {M}achine},
  howpublished = {arXiv preprint},
  year         = {2024},
  eprint       = {2403.09729},
  archivePrefix= {arXiv},
  primaryClass = {math.CA},
  doi          = {10.48550/arXiv.2403.09729},
}

@article{BradfordEtAlChomp2020,
  author  = {Bradford, Alexander and Day, J. Kain and Hutchinson, Laura and Kaperick, Bryan and Larson, Craig E. and Mills, Matthew and Muncy, David and Van Cleemput, Nico},
  title   = {Automated Conjecturing {II}: Chomp and Reasoned Game Play},
  journal = {Journal of Artificial Intelligence Research},
  volume  = {68},
  pages   = {447--461},
  year    = {2020},
  doi     = {10.1613/jair.1.12188}
}

@misc{SchuergerWarnbergYoung2024Zalpha,
  author       = {Schuerger, Houston and Warnberg, Nathan and Young, Michael},
  title        = {Zero Forcing and Vertex Independence Number on Cubic and Subcubic Graphs},
  year         = {2024},
  eprint       = {2410.21724},
  archivePrefix= {arXiv},
  primaryClass = {math.CO},
  url          = {https://arxiv.org/abs/2410.21724}
}

@misc{DavilaSchuergerSmall2024ClawFreeZeqZplus,
  author       = {Davila, Randy and Schuerger, Houston and Small, Ben},
  title        = {A Characterization of Claw-Free Graphs using Zero Forcing Invariants},
  year         = {2024},
  eprint       = {2412.03463},
  archivePrefix= {arXiv},
  primaryClass = {math.CO},
  url          = {https://arxiv.org/abs/2412.03463},
  note          = {Under revision (Journal of Graph Theory)}
}

@misc{totaldom2007,
  author       = {DeLaVi{\~n}a, Ermelinda and Pepper, Ryan and Waller, William},
  title        = {Some Conjectures of {G}raffiti.pc on Total Domination},
  howpublished = {Preprint (online PDF)},
  year         = {2007},
  url          = {https://www.uhd.edu/documents/academics/sciences/total-dom-2007.pdf}
}

@article{Horgan1993DeathOfProof,
  author    = {John Horgan},
  title     = {The Death of Proof},
  journal   = {Scientific American},
  volume    = {269},
  number    = {4},
  year      = {1993},
  pages     = {92--103},
  month     = {October},
  note      = {Contrasts Penrose’s skepticism with Graffiti’s conjecture-making},
  url       = {https://www.jstor.org/stable/24941638}
}

@book{Penrose1989Emperor,
  author    = {Roger Penrose},
  title     = {The Emperor's New Mind: Concerning Computers, Minds and the Laws of Physics},
  publisher = {Oxford University Press},
  year      = {1989},
  address   = {Oxford, UK},
  isbn      = {978-0198519737}
}

@book{Penrose1994Shadows,
  author    = {Roger Penrose},
  title     = {Shadows of the Mind: A Search for the Missing Science of Consciousness},
  publisher = {Oxford University Press},
  year      = {1994},
  address   = {Oxford, UK},
  isbn      = {978-0195106466}
}

@misc{mitchell2011pulp,
  title        = {PuLP: A Linear Programming Toolkit for Python},
  author       = {Mitchell, Stuart and O'Sullivan, Michael and Dunning, Iain},
  year         = {2011},
  note         = {Optimization Online},
  url          = {https://optimization-online.org/wp-content/uploads/2011/09/3178.pdf}
}

@misc{forrest2005cbc,
  title        = {CBC User Guide},
  author       = {Forrest, John and Lougee-Heimer, Robin},
  year         = {2005},
  publisher    = {INFORMS},
  doi          = {10.1287/educ.1053.0020},
  note         = {COIN-OR Branch-and-Cut (CBC) documentation}
}

@ARTICLE{2020SciPy-NMeth,
  author  = {Virtanen, Pauli and Gommers, Ralf and Oliphant, Travis E. and
            Haberland, Matt and Reddy, Tyler and Cournapeau, David and
            Burovski, Evgeni and Peterson, Pearu and Weckesser, Warren and
            Bright, Jonathan and {van der Walt}, St{\'e}fan J. and
            Brett, Matthew and Wilson, Joshua and Millman, K. Jarrod and
            Mayorov, Nikolay and Nelson, Andrew R. J. and Jones, Eric and
            Kern, Robert and Larson, Eric and Carey, C J and
            Polat, {\.I}lhan and Feng, Yu and Moore, Eric W. and
            {VanderPlas}, Jake and Laxalde, Denis and Perktold, Josef and
            Cimrman, Robert and Henriksen, Ian and Quintero, E. A. and
            Harris, Charles R. and Archibald, Anne M. and
            Ribeiro, Ant{\^o}nio H. and Pedregosa, Fabian and
            {van Mulbregt}, Paul and {SciPy 1.0 Contributors}},
  title   = {{{SciPy} 1.0: Fundamental Algorithms for Scientific
            Computing in Python}},
  journal = {Nature Methods},
  year    = {2020},
  volume  = {17},
  pages   = {261--272},
  doi     = {10.1038/s41592-019-0686-2}
}

@article{huangfu2018parallelizing,
  title   = {Parallelizing the dual revised simplex method},
  author  = {Huangfu, Qiang and Hall, J. A. J.},
  journal = {Mathematical Programming Computation},
  year    = {2018},
  volume  = {10},
  number  = {1},
  pages   = {119--142},
  doi     = {10.1007/s12532-017-0130-5}
}

@misc{txgraffiti_pypi,
  title        = {{txgraffiti} (Python package)},
  author       = {Davila, Randy and Eddy, Jillian.},
  howpublished = {Python Package Index (PyPI)},
  year         = {2025},
  url          = {https://pypi.org/project/txgraffiti/},
  note         = {Accessed 2025-12-30}
}

@article{Barioli2010,
  author  = {Barioli, Francesco and Barrett, Wayne and Fallat, Shaun M. and Hall, H. Tracy and Hogben, Leslie and Shader, Bryan and van den Driessche, P. and van der Holst, Hein},
  title   = {Zero forcing parameters and minimum rank problems},
  journal = {Linear Algebra and its Applications},
  volume  = {433},
  number  = {2},
  pages   = {401--411},
  year    = {2010},
  doi     = {10.1016/j.laa.2010.03.008}
}

\end{document}